\documentclass[preprint,12pt,numafflabel]{elsarticle}
\usepackage{lineno}
\usepackage{pgf,tikz}
\usetikzlibrary{arrows}
\usepackage{tikz}
\usepackage[section]{placeins}
\usepackage{array,multirow,makecell}
\setcellgapes{1pt} \makegapedcells

\usepackage{diagbox}
\usepackage{amsmath}
\usepackage{amsmath,amssymb}
\usepackage{amsthm}
\usepackage[colorlinks=true
linkscolor=red, urlcolor=cyan]{hyperref}

\usepackage{subcaption}
\usepackage{xcolor,colortbl}
\usepackage{graphicx}
\usepackage{here}
\usepackage{mathrsfs}
\usepackage{algpseudocode}
\usepackage{epsfig}
\usepackage{color}
\usepackage{pgf,tikz}
\usetikzlibrary{arrows}
\usepackage{bbold}
\usepackage{lmodern}
\renewenvironment{proof}{{\bfseries Proof.}}{\hfill $\Box$ \\}

\usepackage{amsthm,latexsym}
\usepackage{algorithmicx, algpseudocode, algorithm}
\usepackage{mathtools}
\usepackage{titlesec}
\usepackage{subcaption}
\usepackage{caption}
\usepackage{placeins}
\usepackage{multicol}
\titleformat{\subsection}
{\normalfont\bfseries}{}
{0pt}
{ \thesubsection\quad}
\titleformat{\subsubsection}
{\normalfont\bfseries}{}
{0pt}
{ \thesubsubsection\quad}

\theoremstyle{plain}
\newtheorem{theorem}{Theorem}[section]
\newtheorem{lemma}[theorem]{Lemma}
\newtheorem{proposition}[theorem]{Proposition}
\newtheorem{corollary}[theorem]{Corollary}

\theoremstyle{definition}
\newtheorem{definition}[theorem]{Definition}
\newtheorem{remark}[theorem]{Remark}

\newcommand{\comment}[1]{}

\journal{journal}

\usetikzlibrary{arrows,calc}

\definecolor{wcgreen}{RGB}{220,252,231}
\definecolor{wcred}{RGB}{254,226,226}
\definecolor{headgray}{RGB}{30,41,59}
\definecolor{headfg}{RGB}{241,245,249}

\usepackage{longtable}
\usepackage{booktabs}
\usepackage{xcolor}
\usepackage{colortbl}
\usepackage{array}
\usepackage{amsmath}

\usepackage{bbold}

\begin{document}

\begin{frontmatter}

\title{Cohen--Macaulayness of 4-regular Circulant Graphs $C_n(a,b)$}

\author[1]{Subha Shree R\fnref{fn1}}
\author[1]{Bidwan Chakraborty\corref{cor1}}

\fntext[fn1]{The first author acknowledges support from the University Grants Commission (UGC), Government of India.}
\cortext[cor1]{Corresponding author}

\address[1]{%
Department of Mathematics, National Institute of Technology Puducherry,
Karaikal--609609, India\\
E-mail: ma25d1002@nitpy.ac.in, bidwan.c@nitpy.ac.in
}

\begin{abstract}
In this paper, we investigate the well-coveredness of 4-regular
circulant graphs $C_n(a,b)$, and the
Cohen--Macaulayness of their associated edge ideals. We first
establish that $C_n(a,b)$ is well-covered whenever $n=3a$ or
$n=3b$, and that its independence number is $n/3$ for $n>6$.
For the remaining cases, we obtain sufficient arithmetic
conditions for non-well-coveredness in terms of
$g=\gcd(b,n)$ and $m=n/g$. In particular, we prove that
$C_n(a,b)$ is not well-covered if $m\geq4$ is even and
$g\geq6$, or if $m\geq5$ is odd and $g\geq4$.
We then study the Cohen--Macaulayness of the corresponding
edge ideals. Although the families with $n=3a$ or $n=3b$
are well-covered, we show that their edge ideals are not
Cohen--Macaulay over any field when $n=3a$ with $a\geq4$,
or when $n=3b$ with $b\geq5$. Our approach uses a block
decomposition of the graphs and the construction of faces
with disconnected links, together with Reisner's criterion.
Finally, we provide computational verification for the
remaining small parameter cases, completing the
Cohen--Macaulayness classification within the two
well-covered families.
\end{abstract}

\end{frontmatter}

\renewcommand{\thefootnote}{}

\section{Introduction}

Cohen--Macaulayness is a fundamental property in commutative algebra
and algebraic geometry, reflecting a close relationship between the
depth and Krull dimension of a ring or module. For squarefree monomial
ideals, this property is closely connected to combinatorial structures
through the Stanley--Reisner correspondence
\cite{stanley1996combinatorics,villarreal2018monomial}. In particular,
the Cohen--Macaulayness of the quotient ring associated with the edge
ideal of a graph can be studied through the independence complex of
that graph. Since every Cohen--Macaulay simplicial complex is pure,
Cohen--Macaulayness imposes a necessary combinatorial condition on
the graph: all its maximal independent sets must have the same
cardinality. Graphs satisfying this condition are called
\emph{well-covered graphs}, following Plummer
\cite{plummer1970some,plummer1993well}. Thus, well-coveredness provides
a useful combinatorial condition for investigating the
Cohen--Macaulayness of edge ideals. This connection has motivated
substantial research on well-covered and Cohen--Macaulay graphs in
various graph classes, including trees, chordal graphs, and bipartite
graphs. In the bipartite case, Herzog and Hibi obtained a
combinatorial characterization of Cohen--Macaulayness
\cite{herzog2005distributive}.

More precisely, let $G=(V,E)$ be a finite simple graph with vertex
set $V=\{x_1,\ldots,x_n\}$, and let
$R=K[x_1,\ldots,x_n]$ be a polynomial ring over a field $K$.
The \emph{edge ideal} of $G$ is the squarefree monomial ideal
$
 I(G)=\langle x_ix_j:\{x_i,x_j\}\in E\rangle.
$
The \emph{independence complex} of $G$, denoted by $\Delta(G)$,
is the simplicial complex whose faces are the independent sets of
$G$. Its facets are precisely the maximal independent sets of
$G$. Moreover, the Stanley--Reisner ideal of $\Delta(G)$ is
exactly $I(G)$, and hence its Stanley--Reisner ring is
$K[\Delta(G)]=R/I(G)$. Consequently, the purity of
$\Delta(G)$ is equivalent to the well-coveredness of $G$, while
the Cohen--Macaulayness of $R/I(G)$ is equivalent to the
Cohen--Macaulayness of the associated simplicial complex.

An important class of graphs for investigating these questions
is the class of \emph{circulant graphs}. For
$
 S\subseteq
 \left\{1,2,\ldots,\left\lfloor\frac{n}{2}\right\rfloor\right\},$ the circulant graph $C_n(S)$ has vertex set $\mathbb{Z}_n$,
where two distinct vertices $x$ and $y$ are adjacent if
$
 y-x\equiv s\ \text{or}\ -s\pmod n$
 for some $ s\in S.$ Equivalently, circulant graphs are Cayley graphs over cyclic groups. They include cycles and complete graphs as particular
examples. Their regular structure and arithmetic description
make them a natural family for studying the relationship
between graph properties and algebraic properties of edge
ideals. The independence numbers and maximal independent
sets of circulant graphs have recently been investigated by
Carotenuto, Nicolosa, and Salvatore using an
odd-cycle-transversal approach \cite{carotenuto2026maximum},
providing a complementary perspective to the well-coveredness
questions considered here.

 A circulant graph is called \emph{4-regular} if every vertex has exactly four neighbors. In this paper, we focus on 4-regular circulant graphs of the form $G=C_n(a,b)$, where the vertex set is $\mathbb{Z}_n$ and two distinct vertices $i,j\in\mathbb{Z}_n$ are adjacent whenever $i-j\equiv \pm a$ or $\pm b\pmod n$. Throughout the paper, we assume that
$1\leq a<b\leq \left\lfloor\frac{n-1}{2}\right\rfloor$
and $\gcd(n,a,b)=1$. These assumptions ensure that the four residues $\pm a,\pm b$ are distinct and nonzero modulo $n$, so that every vertex has exactly four neighbours. Our objective is to investigate the well-coveredness of these graphs and the Cohen--Macaulayness of their associated edge ideals.

The well-coveredness problem for circulant graphs has been
studied from a computational and combinatorial perspective.
In particular, \cite{BROWN2011244} investigated the complexity
of recognizing well-covered circulant graphs and characterized
several families using independence-polynomial methods. The Cohen--Macaulayness of circulant graphs was studied in
\cite{Vander}, where results were established for several
classes, including $C_n(1,2,\ldots,d)$ and certain cubic
circulant graphs. The independence complexes of well-covered
circulant graphs were further investigated in \cite{Earl2016},
which provided Cohen--Macaulayness characterizations for
several families, including circulants with omitted generators.
More recently, \cite{hoang2021cohen} characterized the
well-coveredness, Cohen--Macaulayness, Buchsbaumness, and
Gorenstein properties of two additional classes of circulant
graphs and their complements.

In this paper, we study the well-coveredness and
Cohen--Macaulayness of $C_n(a,b)$ under the assumption
$\gcd(n,a,b)=1$. Our first main result establishes the
well-coveredness of the two families $n=3a$ and $n=3b$ and
determines their independence number. In particular, for
$n>6$, we obtain $\alpha(C_n(a,b))=n/3$ in both cases.
For the remaining cases, we establish sufficient arithmetic
conditions for non-well-coveredness in terms of
$g=\gcd(b,n)$ and $m=n/g$. More precisely, we prove that
$C_n(a,b)$ is not well-covered if either $m\geq4$ is even
and $g\geq6$, or $m\geq5$ is odd and $g\geq4$.

\vskip 2mm
\noindent\textbf{Theorem \ref{wellcovered}.}
Let $n,a,b\in\mathbb{Z}_{>0}$ satisfy
$1\leq a<b\leq \left\lfloor\frac{n-1}{2}\right\rfloor$
and $\gcd(n,a,b)=1$ with $n\geq19$, and let $G=C_n(a,b)$. Set $g=\gcd(b,n)$ and $m=n/g$. Then:
\begin{enumerate}
\item[(a)] If $n=3a$ or $n=3b$, then $G$ is well-covered and $\alpha(G)=n/3$ for $n>6$.

\item[(b)] If $n\neq 3a$ and $n\neq 3b$, then $G$ is not well-covered whenever
\[
\begin{cases}
m\geq 4\text{ is even}, & g\geq 6,\\
m\geq 5\text{ is odd}, & g\geq 4.
\end{cases}
\]
\end{enumerate}

We then turn to the Cohen--Macaulayness of the associated edge ideals. Although well-coveredness is a necessary condition for Cohen--Macaulayness, it is not sufficient. In particular, we show that the edge ideals of the well-covered families above are not Cohen--Macaulay over any field when $n=3a$ with $a\geq 4$, or when $n=3b$ with $b\geq 5$.

In particular, we give a complete classification for Cohen--Macaulayness when $n=3a$ or $n=3b$ , which shows that the well-covered families $n=3a$ and $n=3b$ are never Cohen--Macaulay over any field.

\vskip 2mm
\noindent\textbf{Theorem \ref{combined}.}
Let $n,a,b\in\mathbb{Z}_{>0}$ satisfy
$1\leq a<b\leq \left\lfloor\frac{n-1}{2}\right\rfloor$
and $\gcd(n,a,b)=1$ with $n\geq19$, and let $G=C_n(a,b)$. Set $g=\gcd(b,n)$ and $m=n/g$. Then:
\begin{enumerate}
\item[(a)] If $n=3a$ or $n=3b$, then $G$ is well-covered and $\alpha(G)=n/3$ for $n>6$. Moreover, if $n=3a$ with $a\geq 4$, or if $n=3b$ with $b\geq 5$, then $R/I(G)$ is not Cohen--Macaulay over any field.

\item[(b)] Suppose that $n\neq 3a$ and $n\neq 3b$. If either
\[
\begin{cases}
m\geq 4\text{ is even}, & g\geq 6,\\
m\geq 5\text{ is odd}, & g\geq 4,
\end{cases}
\]
then $G$ is not well-covered. Consequently, $R/I(G)$ is not Cohen--Macaulay over any field.
\end{enumerate}

Both Theorem~\ref{wellcovered} and Theorem~\ref{combined} are stated under the standing assumption $\gcd(n,a,b)=1$, which is exactly the condition for $C_n(a,b)$ to be connected. This assumption is not restrictive. Theorem~\ref{main} shows that, for $d=\gcd(n,a,b)>1$, the graph $C_n(a,b)$ is a disjoint union of $d$ copies of the connected circulant graph $H=C_{n/d}(a/d,b/d)$, and that well-coveredness and Cohen--Macaulayness of $C_n(a,b)$ are equivalent to the corresponding properties of $H$. Consequently, every result obtained for $\gcd(n,a,b)=1$ transfers to an arbitrary circulant graph $C_n(a,b)$, with the hypotheses imposed on the reduced parameters $(n/d,a/d,b/d)$.

\vskip 2mm
\noindent\textbf{Theorem \ref{main}.}

    Let $d=gcd(n,a,b)$ , $G=C_n(a,b)$ and $H = C_\frac{n}{d} (\frac{a}{d}, \frac{b}{d})$ then,
    \begin{itemize}
        \item[(a)]  $G$ is connected if and only if $d=1$.
        \item[(b)] If $d>1$, then $G\cong d H $ ,  $G$ is isomorphic to $d$ disjoint copies of $C_\frac{n}{d} (\frac{a}{d}, \frac{b}{d})$.
        \item[(c)]  $G$ is well-covered if and only if  $H$ is well covered.
        \item[(d)] $G$ is Cohen--Macaulay if and only if $H$ is Cohen--Macaulay.
        \end{itemize}

The paper is organized as follows. Section~2 recalls
the necessary background on graphs, simplicial
complexes, independence complexes, Stanley--Reisner
ideals, and the relationship between
Cohen--Macaulayness and well-coveredness. Section~3
develops the block decomposition and establishes the
well-coveredness and non-well-coveredness results.
Section~4 uses the well-covered cases to construct
faces whose links are disconnected and applies
Reisner's criterion to prove non-Cohen--Macaulayness.
Section~5 presents computational verification for
the small parameter cases considered in the paper and then it extends the
assumption $\gcd(n,a,b)=1$, to $\gcd(n,a,b)>1$ .

\section{\centering PRELIMINARIES}

We follow standard texts \cite{monomial,bondy1976graph}in these fields.  Throughout the paper, $K$ will denote any arbitrary field, and $R= K[x_1,...,x_n]$ will be a polynomial ring over $K$. For obvious reasons, we shall identify the variables $x_1,...,x_n$ with the vertices of simplicial complexes and the vertices of the graphs being discussed. For a graph $G$, let $V(G)$ denote the set of vertices of $G$, and let
$E(G)$ denote the set of edges of $G$. An edge $e\in E(G)$ connecting two vertices $x$ and $y$ will also be written as $\{x,y\}$.

\begin{definition}
For $ S\subseteq \left\{1,2,\ldots,\left\lfloor\frac{n}{2}\right\rfloor\right\},$ the circulant graph $C_n(S)$ has vertex set $\mathbb{Z}_n$, where two distinct vertices $x$ and $y$ are adjacent if
$ y-x\equiv s\ \text{or}\ -s\pmod n$
 for some $ s\in S.$ Equivalently, $E(C_n(S))=\{\{x,y\}:x,y\in\mathbb{Z}_n,\ y-x\equiv s\ \text{or}\ -s\pmod n$
 for some $ s\in S$ \} .

Particularly in this paper we use the circulant graph $C_n(a,b)$, let $1\leq a<b$ and let $a,b\in\mathbb{Z}_n\setminus\{0\}$. The circulant graph $C_n(a,b)$ is the graph $C_n(\{a,b\})$. Thus, $V(C_n(a,b))=\mathbb{Z}_n$, and two distinct vertices $x,y\in\mathbb{Z}_n$ are adjacent if and only if
$y-x\equiv a$ or $-a$ or $b$ or $-b$ $\pmod n$ .
\end{definition}

\begin{definition}
An edge $\{x,y\}$ of $C_n(a,b)$ is called an \emph{$a$-edge} if $y-x \equiv a$ or $-a$ $\pmod n$ , and it is called a \emph{$b$-edge} if $y-x \equiv b$ or $-b$ $\pmod n$ . Thus, the edges of $C_n(a,b)$ are precisely the $a$-edges and the $b$-edges.
\end{definition}

\begin{definition}
  The \emph{neighborhood} of a vertex $x$ in $G$ is the set
$N_G(x)=\{y\in V(G):\{x,y\}\in E(G)\}$, and the \emph{closed neighborhood} of
$x$ in $G$ is the set $N_G[x]=N_G(x)\cup\{x\}$. The number $\deg_G(x)=|N_G(x)|$
is called the degree of $x$ in $G$.
\end{definition}

\begin{definition}
A subset $S$ of $V$ is called an \emph{independent set} of $G$ if no two vertices of $S$ are adjacent in $G$.  An independent set $S$ is called \emph{maximum} if $G$ has no independent set $S'$ with $|S|>|S'|$. The size of the largest independent set in $G$, denoted by
$\alpha(G)$, is called the independence number of $G$.
\end{definition}

\begin{definition}
    An independent set $S$ is \emph{maximal} if $S \cup \{v\} $ is not independent for any $v \in V \setminus S$
\end{definition}

\begin{definition}
If all maximal independent sets in $G$ have the same size, then $G$ is called  \emph{well-covered}.
\end{definition}

\begin{definition}
    A subset $M\subseteq E$ is a \emph{matching} if every element of $M$ is an edge of $G$, and no two edges in $M$ are adjacent. A vertex $v$ is said to be \emph{$M$-saturated} if some edge of $M$ is incident with $v$. A matching $M$ is called a \emph{perfect matching} if every vertex of $G$ is $M$-saturated. In other words, a perfect matching matches every vertex exactly once.
\end{definition}

\begin{definition}
An \emph{simplicial complex}  $\Delta$ on the vertex set $\{1,\ldots,n\}$ is a collection of subsets called faces or simplices, closed under taking subsets; that is, if $\Lambda\in\Delta$ is a face and $\tau\subseteq\Lambda$, then $\tau\in\Delta$. A simplex
$\Lambda\in\Delta$ of cardinality $|\Lambda|=i+1$ has dimension $i$ and is called an $i$-face of $\Delta$. The dimension $\dim(\Delta)$ of $\Delta$ is the maximum of the dimensions of its faces, or it is $-\infty$ if $\Delta=\{\}$ is the void complex, which has no faces.
\end{definition}

\begin{definition}
The \emph{independence complex} of $G$, denoted by $\operatorname{Ind}(G)$, is the simplicial complex whose faces are precisely the independent sets of
$G$. Thus, $F\in\operatorname{Ind}(G)$ if and only if $F$ is an independent set of $G$.
\end{definition}

\begin{definition}
Let $\Delta$ be a simplicial complex on the vertex set
$X=\{x_1,\ldots,x_n\}$, and let $R=K[x_1,\ldots,x_n]$. The
\emph{Stanley--Reisner ideal} of $\Delta$ is the squarefree monomial
ideal $I_\Delta=\left\langle x^F: F\subseteq X,\ F\notin\Delta\right\rangle
\subseteq R$, generated by the monomials corresponding to the \emph{minimal nonfaces} of $\Delta$, that is, the subsets $F\notin\Delta$ that are minimal with respect to inclusion. The quotient $R/I_\Delta$ is called the \emph{Stanley--Reisner ring} of $\Delta$.

This construction gives a bijective correspondence between simplicial
complexes on $X$ and squarefree monomial ideals in $R$. Under this
correspondence, for a graph $G$, the edge ideal $I(G)$ coincides with the Stanley--Reisner ideal $I_{\operatorname{Ind}(G)}$ of its independence complex.
\end{definition}

\begin{definition}
Let $R=K[x_1,\ldots,x_n]$ with maximal homogeneous ideal
$\mathfrak{m}=(x_1,\ldots,x_n)$, and let $M$ be a finitely generated
graded $R$-module. The \emph{depth} of $M$, denoted by
$\operatorname{depth}(M)$, is the length of a maximal $M$-regular
sequence of elements of $\mathfrak{m}$. Since
$\operatorname{depth}(M)\leq\dim(M)$, the module $M$ is called
\emph{Cohen--Macaulay} if $\operatorname{depth}(M)=\dim(M)$. For a simplicial complex $\Delta$, we say that $\Delta$ is
\emph{Cohen--Macaulay over $K$} if the Stanley--Reisner ring
$R/I_\Delta$ is Cohen--Macaulay as an $R$-module. A graph $G$ is called \emph{Cohen--Macaulay over $K$} if its edge ideal $I(G)$, equivalently $R/I(G)$, has this property.
\end{definition}

Let $\Delta$ be a simplicial complex on $\{1,\ldots,n\}$. For each
integer $i$, let $F_i(\Delta)$ be the set of $i$-faces of $\Delta$. For each $i$, let $F_i(\Delta)$ be the set of $i$-faces of $\Delta$. Let $K^{F_i(\Delta)}$ be a vector space over $K$ .

 \begin{definition}
The reduced chain complex of $\Delta$ over $k$ is the
complex  $C_{\bullet}(\Delta;K)$:
$0\leftarrow K^{F_{-1}(\Delta)}
\xleftarrow{\partial_0}\cdots
\xleftarrow{\partial_i} K^{F_{i}(\Delta)}
\xleftarrow{\partial_{i+1}}\cdots
\xleftarrow{\partial_{n-1}} K^{F_{n-1}(\Delta)}\leftarrow0$. For each integer $i$, the vector space
$\widetilde{H}_i(\Delta;K)=\ker(\partial_i)/\operatorname{im}(\partial_{i+1})$ is called the \emph{$i$th reduced homology} of $\Delta$ over $K$.

The dimension  $\widetilde{H}_0(\Delta;K)$, as a $K$-vector space is one less than the number of connected components of $\Delta$. Elements of $\ker(\partial_i)$ are often called \emph{$i$-cycles} and elements of $\operatorname{im}(\partial_{i+1})$ are often called
\emph{$i$-boundaries}.

\end{definition}

\begin{definition}
Let $\Delta$ be a simplicial complex and let $F\in\Delta$. The \emph{link} of $F$ in $\Delta$ is $\operatorname{lk}_{\Delta}(F) =\{H\in\Delta:H\cap F=\emptyset,\ F\cup H\in\Delta\}$.
\end{definition}

\begin{definition}
 Let $X$ and $Y$ be simplicial complexes on disjoint vertex sets. The \emph{join} of $X$ and $Y$ is the simplicial complex
$X*Y=\{\sigma\cup\tau \;:\; \sigma\in X,\ \tau\in Y\}.$ Each face of $X*Y$ can be written uniquely as $\sigma\cup\tau$ with $\sigma\in X$ and $\tau\in Y$, and
$\dim(\sigma\cup\tau)=\dim\sigma+\dim\tau+1$. In particular,
$\dim(X*Y)=\dim X+\dim Y+1$
\end{definition}

\section{\centering WELL COVEREDNESS OF $C_n(a,b)$}

In this section, we determine the well-coveredness of $C_n(a,b)$. Let $G=C_n(a,b)$, where
$1\leq a<b\leq\left\lfloor\frac{n-1}{2}\right\rfloor$, with vertex set
$V(G)=\mathbb{Z}_n$. For $c\in\{a,b\}$, define
$g_c=\gcd(c,n)$ and $m_c=n/g_c$. Throughout the paper, we use
$g=g_b=\gcd(b,n)$ and $m=m_b=n/g$. Put
$t=b/g$. Then $b=gt$ and
$\gcd(t,m)=1$, so $t$ is invertible modulo $m$.

We partition the vertex set $V(G)=\mathbb{Z}_n$ into the blocks
$B_r=\{r+kg:0\leq k<m\}$, where
$r=0,1,\ldots,g-1$. These blocks are pairwise disjoint; their union is
$\mathbb{Z}_n$, and each block contains exactly $m$ vertices. The
$b$-edges preserve each block and form a copy of the cycle $C_m$ within each block, while the $a$-edges connect different blocks. We use this decomposition to identify the
cases in which $C_n(a,b)$ is well-covered, namely the exceptional cases $n=3a$ and $n=3b$. Finally, we construct explicit admissible cyclic sequences to establish non-well-coveredness for a broad range of the remaining parameters

\begin{lemma}\label{Br blocks are cycles}
The $b$-edges induced on each $B_r$ form a cycle $C_m$.
\end{lemma}

\begin{proof}
Fix $r\in\{0,1,\ldots,g-1\}$ and let $r+kg\in B_r$. Since $b=gt$, we have
$r+kg+b=r+(k+t)g\in B_r$
in $\mathbb{Z}_n$. Similarly, $r+kg-b\in B_r$. Thus, every $b$-edge has both endpoints in the same block. Starting from $r+kg$ and repeatedly following $b$-edges, the corresponding indices in $\mathbb{Z}_m$ are
$k,\ k+t,\ k+2t,\ldots,k+(m-1)t$.
Since $\gcd(t,m)=1$, these indices are pairwise distinct modulo $m$. Hence, every vertex of $B_r$ is visited exactly once, and the $b$-edges induced on $B_r$ form a cycle $C_m$. Moreover, under the relabeling
$k\longmapsto t^{-1}k\pmod m$,
a $b$-edge $k\mapsto k+t$ becomes
$t^{-1}k\longmapsto t^{-1}(k+t)=t^{-1}k+1\pmod m$.
Thus, each $B_r$ is identified with the standard cycle $C_m$.
\end{proof}

\begin{lemma}\label{location of a edges}
 Let $B_r$ be the block
decomposition of the vertex set. Then every $a$-edge joins a vertex of $B_r$ to a vertex of $B_{r+r_a\pmod g}$.
\end{lemma}

\begin{proof}
Fix $r\in\{0,\ldots,g-1\}$ and let $r+kg\in B_r$. Since
$a=qg+r_a$, we have $r+kg+a=(r+r_a)+(k+q)g$. If $r+r_a<g$,  the resulting vertex lies in $B_{r+r_a}$. If $r+r_a\geq g$, then
$r+r_a=(r+r_a-g)+g$, and therefore
$r+kg+a=(r+r_a-g)+(k+q+1)g$,
which lies in $B_{r+r_a-g}$. Thus, in both cases, the resulting block index is
$r+r_a\pmod g$. Since $r+kg\in B_r$ was arbitrary, adding $a$ sends every vertex of $B_r$ to a vertex of $B_{r+r_a\pmod g}$. Consequently, every $a$-edge joins two blocks whose indices differ by $r_a$ modulo $g$.
\end{proof}

\begin{lemma}\label{Cyclic order of the blocks}
 The $a$-edges connect the $g$ blocks in the cyclic order $B_0,B_{r_a},B_{2r_a},\ldots,B_{(g-1)r_a}$, where the indices are taken modulo $g$, and each block is visited exactly once.
\end{lemma}

\begin{proof}
Since $g / b$ and $\gcd(n,a,b)=1,$ we have $\gcd(a,g)=1$. Since $r_a\equiv a\pmod g$, it follows that
$\gcd(r_a,g)=1$. Thus $r_a$ is invertible modulo $g$, and the residues $0,r_a,2r_a,\ldots,(g-1)r_a$ are pairwise distinct modulo $g$. Consequently, these residues constitute all the elements of
$\mathbb{Z}_g$. By Lemma~\ref{location of a edges}, every $a$-edge joins a vertex of
$B_j$ to $B_{j+r_a\pmod g}$. Starting from $B_0$, we have the next blocks as
$B_0,B_{r_a},B_{2r_a},\ldots,B_{(g-1)r_a}$.
Since the corresponding residues are pairwise distinct modulo $g$, all $g$ blocks are visited exactly once. The next transition returns to $B_0$, as ${gr_a} \equiv0 \pmod{g}$.
\end{proof}

\begin{lemma} \label{shift}

For each $r$, let $s_r\in
\mathbb{Z}_m$ denote the shift induced by an $a$-edge from $B_r$ to $B_{r+r_a\pmod g}$. Then
$s_r=
\begin{cases}
\ell, & r+r_a<g,\\
\ell+t^{-1}, & r+r_a\geq g,
\end{cases}$
where $t=b/g$ and $\ell\equiv t^{-1}q\pmod m$  .

\end{lemma}

\begin{proof}
From the two cases in the proof of Lemma~\ref{location of a edges}, the second coordinate, that is, the multiple of $g$, increases by $q$ if $r+r_a<g$ and increases by $q+1$ if $r+r_a\geq g$ when the $a$-edge crosses from $B_r$ to $B_{r+r_a\pmod g}$. Looking back at the two cases from Lemma~\ref{location of a edges}, adding $a$ moves a vertex from position $k$ inside its block to position $k+q$ or $k+q+1$, depending only on the block $B_r$ in which the vertex starts,
and not on the particular position $k$. In other words, every vertex in $B_r$ moves by the same amount when it crosses into $B_{r+r_a\pmod g}$. Thus, this movement can be represented by a single
number associated with the entire block. Hence, the induced shift on the second coordinate is constant for all vertices in $B_r$.

Under the relabeling $k\mapsto t^{-1}k\pmod m$, an increment of $q$ becomes
$t^{-1}q\pmod m$, while an increment of $q+1$ becomes
$t^{-1}(q+1)\equiv t^{-1}q+t^{-1}\pmod m$. Writing $\ell\equiv t^{-1}q\pmod m$ gives
$s_r=
\begin{cases}
\ell, & r+r_a<g,\\
\ell+t^{-1}, & r+r_a\geq g.
\end{cases}$ Thus $s_r$ is a well-defined element of  $\mathbb{Z}_m$  associated with the transition $B_r\to B_{r+r_a\pmod g}$ and is independent of the
choice of $k$.
\end{proof}

\begin{theorem}\label{thm:well-covered}
Let $1\leq a<b\leq \left\lfloor\frac{n-1}{2}\right\rfloor$ with
$\gcd(n,a,b)=1$ and $n>6$. If $n=3a$ or $n=3b$, then $C_n(a,b)$ is well-covered. Moreover, $\alpha(C_n(a,b))=\frac{n}{3}$.
\end{theorem}

\begin{proof}
We first consider the case $n=3b$. For $j=0,1,\ldots,b-1$, define $T_j=\{j,j+b,j+2b\}$. Since $n=3b$, the sets $T_0,T_1,\ldots,T_{b-1}$ are pairwise disjoint and form a partition of $V(C_n(a,b))$. Moreover, the difference between any two distinct vertices of $T_j$ is congruent to $b$ or $-b$ modulo $n$. Hence, $T_j$ induces a triangle $K_3$. Therefore, every independent set contains at most one vertex from each
$T_j$, and hence $\alpha(C_n(a,b))\leq b$.

Let $I$ be a maximal independent set. We show that $I\cap T_j\neq\varnothing$ for every $j$. Suppose, to the contrary, that $I\cap T_j=\varnothing$ for some $j$. Since $T_j$ is a triangle, no vertex of $T_j$ can be adjacent to a vertex of $I$ lying in $T_j$. Thus, every vertex of $T_j$ must have a neigbour in $I$, reached by an $a$-edge. The $a$-edges from $T_j$ join $T_j$ to $T_{j+a}$ and $T_{j-a}$, where the subscripts are taken modulo $b$. More precisely, the three edges $j\sim j+a$, $j+b\sim j+a+b$, and $j+2b\sim j+a+2b$ form a perfect matching between $T_j$ and $T_{j+a}$. Similarly, the three $a$-edges between $T_j$ and $T_{j-a}$ form a perfect matching. The two triangles $T_{j+a}$ and $T_{j-a}$ are distinct. Otherwise, $2a\equiv0\pmod b$. Since $\gcd(n,a,b)=1$, this would imply $b/2$, which is impossible because $n=3b > 6$.

Since each of $T_{j+a}$ and $T_{j-a}$ is a triangle, the independent set $I$ contains at most one vertex in each of them. Furthermore, because the $a$-edges form perfect matchings, each selected vertex in these two triangles is adjacent to at most one vertex of $T_j$. Hence, at most two vertices of $T_j$ can be dominated by the vertices of $I$ that are outside $T_j$. Since $T_j$ contains three vertices, at least one vertex $v\in T_j$ has no neighbor in $I$. Because $I\cap T_j=\varnothing$, the set $I\cup\{v\}$ is independent, contradicting the maximality of $I$. Thus, $I\cap T_j\neq\varnothing$ for every $j$. Since $T_j$ is a triangle and $I$ is independent, we have $|I\cap T_j|=1$ for every $j$. Consequently, $|I|=\sum\limits_{j=0}^{b-1}|I\cap T_j|=b$. Hence, every maximal independent set has cardinality $b$.\\

Similarly, we can prove for $n=3a$ case. So that we will be having , $|I|=\sum\limits_{j=0}^{a-1}|I\cap T_j|=a$ for this case. Thus, every maximal independent set has cardinality $a$.  In either case, every maximal independent set has cardinality $n/3$.
Hence $C_n(a,b)$ is well-covered.
\end{proof}

We now investigate the converse direction that is whether well-coveredness of $C_n(a,b)$ forces $n=3a$ or $n=3b$. Thus, suppose that $n\neq 3a$ and $n\neq 3b$. We aim to construct two maximal independent sets of $C_n(a,b)$ having different cardinalities.

Let $B_0,B_{r_a},B_{2r_a},\ldots,B_{(g-1)r_a}$, where the indices are taken modulo $g$. Accordingly, for each $j\in\mathbb Z_g$, define $C_j=B_{jr_a\pmod g}$ and $\sigma_j=s_{jr_a\pmod g}$, where $\sigma_j$ denotes the shift associated with the transition from $C_j$ to $C_{j+1}$. For each $j$, let $S_j$ be an independent set of the copy of $C_m$ induced by $C_j$, expressed in its local coordinates in $\mathbb Z_m$.

\begin{definition}
 For an independent set $S\subseteq V(G)$, its closed neighborhood is $N[S]=\bigcup_{v\in S}N[v]$,
where $N[v]$ denotes the closed neighborhood of the vertex $v$. The \emph{defect set} of $S$ is $D(S)=V(G)\setminus N[S]$. Thus, $D(S)$ consists precisely of the vertices that are not dominated by $S$.
\end{definition}

\begin{definition}\label{def:admissible-sequence}
A sequence $(S_0,S_1,\ldots,S_{g-1})$ of independent sets of $C_m$ is an admissible cyclic sequence if, for every $j\in\mathbb Z_g$,
\begin{enumerate}
\item[(i)] $S_j\cap(S_{j-1}+\sigma_{j-1})=\varnothing$;
\item[(ii)] $S_j\cap(S_{j+1}-\sigma_j)=\varnothing$;
\item[(iii)] $D(S_j)\subseteq(S_{j-1}+\sigma_{j-1})\cup(S_{j+1}-\sigma_j)$.
\end{enumerate}
Here, all block indices are taken modulo $g$, and the translations are taken in $\mathbb Z_m$.
\end{definition}

\begin{remark}
Conditions (i) and (ii) ensure that the selected vertices in consecutive blocks do not contain an edge between them. Condition (iii) ensures that every vertex not selected in a block is dominated by a selected vertex in one of the two neighboring blocks. Consequently, an admissible cyclic sequence determines a maximal independent set of $C_n(a,b)$. Therefore, to prove that $C_n(a,b)$ is not well-covered, it suffices to construct two admissible cyclic sequences whose total cardinalities are different.
\end{remark}

\begin{lemma}\label{lem:mc3}
For $c\in\{a,b\}$, let $g_c=\gcd(c,n)$ and $m_c=\frac{n}{g_c}$. Assume $1\leq a<b\leq \left\lfloor\frac{n-1}{2}\right\rfloor$. Then $m_c=3$ if and only if  $n=3c$.
\end{lemma}

\begin{proof}
Suppose first that $m_c=3$. Then $n=3g_c$. Since $g_c=\gcd(c,n)$, we may write
$c=g_ck$ for some positive integer $k$. Since $c\leq\frac{n-1}{2}$, we have $g_ck\leq\frac{3g_c-1}{2}$. Thus, $2g_ck\leq3g_c-1$, and hence $g_c(2k-3)\leq-1$. As $g_c>0$, this implies $2k-3<0$, so $k<\frac{3}{2}$. Since $k$ is a positive integer, $k=1$. Therefore $c=g_c$. Since $n=3g_c$, we obtain $n=3c$.

Conversely, suppose that $n=3c$. Then $g_c=\gcd(c,n)=\gcd(c,3c)=c$, and consequently $m_c=\frac{n}{g_c}=\frac{3c}{c}=3$. Therefore, $m_c=3$ if and only if $n=3c$.
\end{proof}

By the preceding result, it remains to consider the case $m_c\neq 3$. We shall show that whenever $m_c\neq 3$, the graph has two maximal independent sets of different cardinalities. Hence $C_n(a,b)$ is not well-covered. The only remaining case is $m_c=3$, which, by the lemma~\ref{lem:mc3}, is equivalent to $n=3c$, and it is done in the theorem~\ref{thm:well-covered}. However, from now on, we use only $g_b=g$, so we will use the symbol $g$ throughout the paper. Unless otherwise stated, $g$ refers to $g_b$. Therefore, we first prove the result for the case $n=3b$. For $n=3a$, we interchange the roles of $a$ and $b$. First, we note that $m\geq3$. Indeed, if $m=1$, then
$g=n$, so $n\mid b$, which is impossible since $0<b<n$. If $m=2$, then $g=n/2$, so $(n/2)\mid b$. Since $b>0$, this implies $b\geq n/2$, contradicting $b\leq\left\lfloor\frac{n-1}{2}\right\rfloor<n/2$. Hence $m\geq3$.

\begin{theorem}\label{weight preserving}

The map $I\longmapsto(I\cap C_0,I\cap C_1,\ldots,I\cap C_{g-1})$, with each intersection expressed in the local coordinates of its block, defines a weight-preserving bijection between the maximal independent sets of $G$ and the admissible cyclic sequences. Moreover, if $I\longmapsto(S_0,S_1,\ldots,S_{g-1})$, then
$|I|=\sum\limits_{j=0}^{g-1}|S_j|$.
\end{theorem}

\begin{proof}
Let $I$ be a maximal independent set of $G$. For each $j\in\mathbb Z_g$, let $S_j$ be the  representation of $I\cap C_j$. Since each $C_j$ induces a copy of $C_m$, every $S_j$ is an independent set of $C_m$. The $a$-edges joining $C_{j-1}$ to $C_j$ are described by the shift $\sigma_{j-1}$. Hence, the independence of $I$ implies $S_j\cap(S_{j-1}+\sigma_{j-1})=\varnothing$.
Similarly, the $a$-edges joining $C_j$ to $C_{j+1}$ give
$S_j\cap(S_{j+1}-\sigma_j)=\varnothing$.
Thus, conditions (i) and (ii) of Definition~\ref{def:admissible-sequence} hold. Now let $x\in D(S_j)$. The vertex represented by $x$ is not dominated by $S_j$ within its own block. Since $I$ is maximal, every vertex outside $I$ has a neighbor in $I$. Therefore, the vertex represented by $x$ must be adjacent to a selected vertex in either $C_{j-1}$ or $C_{j+1}$. This gives
$x\in(S_{j-1}+\sigma_{j-1})\cup(S_{j+1}-\sigma_j)$.
Consequently, $D(S_j)\subseteq(S_{j-1}+\sigma_{j-1})\cup(S_{j+1}-\sigma_j)$. Thus, condition (iii) holds, and the resulting sequence is admissible.

Conversely, let $(S_0,S_1,\ldots,S_{g-1})$ be an admissible cyclic sequence. Interpret each $S_j$ as a set of vertices in $C_j$, and let
$I=\bigcup\limits_{j=0}^{g-1}S_j$.
Each $S_j$ is independent within its block, and conditions (i) and (ii) ensure that no $a$-edge joins two selected vertices in consecutive blocks. Hence, $I$ is independent in $G$. Let $v\notin I$ lie in block $C_j$. If $v$ is dominated by $S_j$ within its block, then it has a neighbor in $I$. Otherwise, $v$ belongs to $D(S_j)$ for some block $C_j$. By condition (iii),
$D(S_j)\subseteq(S_{j-1}+\sigma_{j-1})\cup(S_{j+1}-\sigma_j)$, so $v$ has a neighbor in one of the neighboring blocks. Thus, every vertex outside $I$ has a neighbor in $I$, proving that $I$ is maximal. The two constructions are inverse to each other, so they define a bijection. Since the blocks are pairwise disjoint, we have
$|I|=\sum\limits_{j=0}^{g-1}|S_j|$. Hence, the bijection is weight-preserving.
\end{proof}

\subsection{Non well-coveredness when $m$ is even}

Having reduced well-coveredness to a question about admissible cyclic
sequences, we now consider the case where the block length $m$ is even.
We use the two alternating maximum independent sets of $C_m$, namely
$S$ and $S'$. Whether each shift $\sigma_j$ is even or odd determines which
of these two sets can be placed in each block.
Assume $m$ is even. Our aim is to create two admissible sequences with different cardinalities. We use two alternating maximal independent sets of $C_m$ namely $S$ and $S'$ . Let $S=\{0,2,\ldots,m-2\}$ and $S'=\{1,3,\ldots,m-1\}$. Since $m$ is even,
$|S|=|S'|=m/2$ and $D(S)=D(S')=\varnothing$. Thus, if $S_j=S$ or $S_j=S'$, then maximality condition of the admissible cyclic sequence is
automatically satisfied. Write $S^0=S$, $S^1=S'$. Since $m$ is even, the parity $\pi(s)\in\{0,1\}$ of $s\in\mathbb{Z_m}$ is well defined, and for all $x\in\{0,1\}$ and $s\in \mathbb{Z_m}$
\begin{equation}\label{eq:translate}
S^x+s=S^{x+\pi(s)},S^x-s=S^{x+\pi(s)} (\text{indices mod }2).
\end{equation}

Indeed, an even translation fixes each of $S,S'$ and an odd translation interchanges them.

\begin{lemma}
Assume that $m$ is even and $g>1$. Then among the shifts
$\sigma_0,\sigma_1,\ldots,\sigma_{g-1}$, there is at least
one even shift and at least one odd shift.
\end{lemma}
\begin{proof}
Since $g>1$ and $\gcd(r_a,g)=1$, we have
$0<r_a<g$. By the construction of the shifts, exactly
$r_a$ of the shifts are congruent to $\ell+t^{-1}\pmod m$,
while the remaining $g-r_a$ shifts are congruent to
$\ell\pmod m$. Since $m$ is even and $\gcd(t,m)=1$, the integer $t$ is odd, and its multiplicative inverse modulo $m$ is also odd. Consequently, $\ell$ and $\ell+t^{-1}$ have opposite parities. Since $0<r_a<g$, both shift values occur. Hence, among $\sigma_0,\ldots,\sigma_{g-1}$, there is at least one even shift and at least one odd shift.
   \end{proof}

\begin{remark}\label{o,e}
    For each $j$, define $\pi_j= 0,\text{if }\sigma_j\text{ is even},  \pi_j =1,\text{if }\sigma_j\text{ is odd}.
$  Let $o=\sum\limits_{j=0}^{g-1}\pi_j$ be the number of odd shifts. Since there are $g$ shifts altogether, the number of even shifts is $e=g-o$. Thus, $o+e=g$, $1\le o,e\le g-1$.
\end{remark}
\begin{lemma}\label{even shift}
Let $S_j=S^{x_j}$ and $S_{j+1}=S^{x_{j+1}}$, where $x_j,x_{j+1}\in\{0,1\}$. Then the independence conditions between $C_j$ and $C_{j+1}$ hold if and only if $x_{j+1}\equiv x_j+\pi_j+1\pmod 2.$
\end{lemma}

\begin{proof}
The conditions between $C_j$ and $C_{j+1}$ are
$S_{j+1}\cap(S_j+\sigma_j)=\emptyset$ (condition~(i) at $j+1$) and $S_j\cap(S_{j+1}-\sigma_j)=\emptyset$ (condition~(ii) at $j$)  the two are equivalent. By \eqref{eq:translate},
$S_j+\sigma_j=S^{x_j+\pi_j}$. Two sets among $S^0,S^1$ are disjoint if and only if they are different, so the condition holds if and only if $x_{j+1}\not\equiv x_j+\pi_j \pmod{2},$ that is, $x_{j+1}\equiv x_j+\pi_j+1\pmod{2}$. We check the two
possible cases  of $\sigma_j$.

\emph{Case 1: $\sigma_j$ is even.}
Then $\pi_j=0$, and translation by $\sigma_j$, or by $-\sigma_j$, preserves the two sets, $S^0+\sigma_j=S^0$, $S^1+\sigma_j=S^1$, and likewise $S^0-\sigma_j=S^0$, $S^1-\sigma_j=S^1$. Hence two neighboring blocks carry different sets, i.e.\ $x_{j+1}\equiv x_j+1\pmod{2}$. Since $\pi_j=0$, this is precisely $x_{j+1}\equiv x_j+\pi_j+1\pmod{2}$.

\emph{Case 2: $\sigma_j$ is odd.} Then $\pi_j=1$, and translation by $\sigma_j$, or by $-\sigma_j$, interchanges the two sets, $S^0+\sigma_j=S^1$, $S^1+\sigma_j=S^0$, and likewise for $-\sigma_j$. Hence two neighboring blocks carry the same set, i.e.\ $x_{j+1}\equiv x_j\pmod{2}$. Since $\pi_j=1$, $x_j+\pi_j+1\equiv x_j+1+1\equiv x_j\pmod{2}$, so again $x_{j+1}\equiv x_j+\pi_j+1\pmod{2}$.

In both cases, the condition is $x_{j+1}\equiv x_j+\pi_j+1\pmod{2}$. This is exactly the relation that allows two non-empty neighboring blocks to occur in an admissible cyclic sequence.
\end{proof}

\begin{lemma}\label{even admissible}
Let $S_{j-1}=S^{x_{j-1}}$, $S_{j+1}=S^{x_{j+1}}$ and $S_j=\emptyset$, where $x_{j-1},x_{j+1}\in\{0,1\}$. Then, $C_j$ can be inserted in between two non empty blocks in an admissible cyclic sequence  if and only if $x_{j-1}+x_{j+1}\equiv 1+\pi_{j-1}+\pi_j \pmod{2}.$
\end{lemma}

\begin{proof}
Since $S_j=\emptyset$, each of conditions~(i) and~(ii) at the pairs $(C_{j-1},C_j)$ and $(C_j,C_{j+1})$ reads $\emptyset\cap(\cdots)=\emptyset$ and is therefore vacuous. For condition~(iii) at $C_j$, we have $D(S_j)=D(\emptyset)=\mathbb{Z_m}$, and the
right-hand side $(S_{j-1}+\sigma_{j-1})\cup(S_{j+1}-\sigma_j)$ is contained in $\mathbb{Z_m}$. Hence condition~(iii) is equivalent to
\begin{equation}\label{eq:empty-union}
  \mathbb{Z_m}=\; (S_{j-1}+\sigma_{j-1})\cup(S_{j+1}-\sigma_j).
\end{equation}
By \ref{eq:translate}, $S_{j-1}+\sigma_{j-1}=S^{\,x_{j-1}+\pi_{j-1}},$ and, since $-\sigma_j$ has the same parity as $\sigma_j$, $\pi(-\sigma_j)=\pi(\sigma_j)=\pi_j$, so $S_{j+1}-\sigma_j=S^{\,x_{j+1}+\pi_j}.$ Thus the right-hand side of \eqref{eq:empty-union} is
$S^{\,x_{j-1}+\pi_{j-1}}\cup S^{\,x_{j+1}+\pi_j}$. Each of these two sets has $m/2$ elements, and $S^0$ and $S^1$ are complementary in $\mathbb{Z_m}$; consequently their union equals $\mathbb{Z_m}$ if and only if the two sets are
distinct, that is, if and only if $x_{j-1}+\pi_{j-1}\;\not\equiv\; x_{j+1}+\pi_j \pmod{2}.$ This is equivalent
to $x_{j-1}+x_{j+1}\;\equiv\; 1+\pi_{j-1}+\pi_j \pmod{2},$ which is the stated congruence.

\end{proof}
\begin{proposition}\label{complete even}
Let $G$ be the circulant graph with an even block size $m$ and $g \ge 3$  blocks $C_0, \dots, C_{g-1}$, with shifts $\sigma_0, \dots, \sigma_{g-1} \in
\mathbb{Z}_m$. Put $\pi_j = \pi(\sigma_j)$, and let $e$ and $o = g - e$ be the number of even and odd shifts, respectively. Let $P = \{k_1, \dots, k_p\} \subseteq \mathbb{Z}_g$ be a set of $p$ pairwise non-adjacent block
indices , with $0 \le p \le \lfloor g/2 \rfloor$. Then there exists a maximal independent set of $G$ whose intersection with block $C_j$ is empty for $j \in P$, and is one of the two alternating maximum independent sets of $C_m$ for $j \notin P$, if and
only if $p \equiv e \pmod 2$. Whenever such a set exists, it has cardinality $(g - p)m/2$.
\end{proposition}

\begin{proof}
For $j\notin P$ write $S_j=S^{x_j}$ with $x_j\in\{0,1\}$. Since
$D(S^0)=D(S^1)=\emptyset$, Lemmas~\ref{even shift}
and~\ref{even admissible} together characterise admissibility of a sequence.

\emph{Case $p=0$.} Then $P$ is empty, so by Lemma~\ref{even shift} the congruence $x_{j+1}\equiv x_j+\pi_j+1\pmod{2}$ must hold for every $j\in\mathbb{Z_m}$. Summing these $g$ congruences,
 $ 0\equiv\sum\limits_{j\in\mathbb{Z_g}}(\pi_j+1)=o+g\equiv e\pmod{2},$ and hence $e$ must be even. Conversely, if $e$ is even, choose
$x_0\in\{0,1\}$ arbitrarily and define $x_1,\dots,x_{g-1}$ successively by the relation of Lemma~\ref{even shift} for $j=0,\dots,g-2$; the last congruence for $j=g-1$ then holds because the sum of all $g$ congruences is consistent.

\emph{Case $p\ge1$.} List the elements of $P$ in cyclic order as
$k_1<\dots<k_p$ and put $k_{p+1}=k_1+g$. By pairwise non-adjacency, $k_{i+1}-k_i\ge2$, so between consecutive elements of $P$ lie the non-empty blocks $k_i+1,\dots,k_{i+1}-1$. Put $z_i=x_{k_i+1}$ and $z_{p+1}=z_1$. Applying Lemma~\ref{even shift} successively from $k_i+1$ to $k_{i+1}-2$,
\begin{equation}\label{eq:run}
  x_{k_{i+1}-1}\equiv z_i+\sum\limits_{j=k_i+1}^{k_{i+1}-2}(\pi_j+1)\pmod{2}.
\end{equation}
At the empty block $k_{i+1}$, Lemma~\ref{even admissible} gives
$x_{k_{i+1}-1}+z_{i+1}\equiv1+\pi_{k_{i+1}-1}+\pi_{k_{i+1}}\pmod{2}$. Substituting \eqref{eq:run} and reducing modulo $2$,
\begin{equation}\tag{$E_i$}
  z_i+z_{i+1}\equiv\beta_i:=1+\sum\limits_{j=k_i+1}^{k_{i+1}}\pi_j+(k_{i+1}-k_i)\pmod{2},\qquad i=1,\dots,p.
\end{equation}
Thus the sequence is admissible if and only if the $z_i$ satisfy
$(E_1),\dots,(E_p)$ and the remaining $x_j$ are obtained by iterating Lemma~\ref{even shift} along each run. Adding $(E_1),\dots,(E_p)$, the left-hand side is $0$, because each $z_i$ occurs exactly twice. On the right-hand side, the integer intervals $\{k_i+1,\dots,k_{i+1}\}$ (with $k_{p+1}=k_1+g$) cover $\mathbb{Z_g}$ exactly once, so
$\sum\limits_i\sum\limits_{j=k_i+1}^{k_{i+1}}\pi_j=\sum\limits_{j\in\mathbb{Z_g}}\pi_j=o$, and
$\sum\limits_i(k_{i+1}-k_i)=g$. Hence $ 0\equiv p+o+g\equiv p+e\pmod{2},$ that is, $p\equiv e\pmod{2}$.

Suppose $p\equiv e\pmod{2}$. Choose $z_1\in\{0,1\}$ arbitrarily. If $p\ge2$, the congruences $(E_1),\dots,(E_{p-1})$ successively determine $z_2,\dots,z_p$; if $p=1$ there is nothing to determine. When we reach the final equation $(E_p)$, it will automatically hold true. To see why, if we add the right-hand sides of all $p$ equations together, we get $p+o+g$. Since $o+g \equiv e$ and our assumption is $p \equiv e$, this sum equals $2p$, which is an even number. Because the total sum is zero, the cycle of equations closes perfectly without any contradictions. Once the $z_i$ are found, we use Lemma~\ref{even shift} to fill in the remaining labels $x_j$ along each run. Then the hypotheses of Lemma~\ref{even shift} and Lemma~\ref{even admissible} are satisfied at every transition, so the sequence is admissible.

Finally, there are $g-p$ non-empty blocks, each carrying
$|S^0|=|S^1|=m/2$ vertices, so the corresponding maximal independent set of $G$ has $(g-p)m/2$ elements.
\end{proof}

\begin{theorem}\label{even case}
Let $m$ be even, let $\sigma_0,\dots,\sigma_{g-1}\in\mathbb{Z_m}$, and let $e$ be the number of even shifts among $\sigma_0,\dots,\sigma_{g-1}$. Then
$G$ is not well-covered if $g\ge4$ and $e$ is even, or if $g\ge6$ and $e$ is odd. In particular, $G $ is not well-covered whenever $g\ge6$. Consequently, $G=C_n(a,b)$ is not well-covered whenever $g_c\ge6$ for some $c\in\{a,b\}$ with $m_c$ even.
\end{theorem}

\begin{proof}
By Proposition~\ref{complete even}, it suffies to exhibit two maximal independent sets of $G$ with different cardinalities.

\emph{Case 1: $e$ even and $g\ge4$.} Since
$e\equiv0\equiv2\pmod{2}$, Proposition~\ref{complete even} with
$P=\emptyset$ (so $p=0$) gives a maximal independent set of size
$W_0=g\,m/2$, and with $P=\{0,2\}$ (so $p=2$) gives an admissible
sequence of size $W_2=(g-2)m/2$. The indices $0$ and $2$ are not
cyclically consecutive because $g\ge4$, and $p=2\le\lfloor g/2\rfloor$,
so $P$ satisfies the hypotheses of Proposition~\ref{complete even}. Applying the proposition gives a maximal independent set of size $W_2=(g-2)m/2$.  Then
$W_0-W_2=m>0$.

\emph{Case 2: $e$ odd and $g\ge6$.} Since
$e\equiv1\equiv3\pmod{2}$, Proposition~\ref{complete even} with
$P=\{0\}$ (so $p=1$) gives a maximal independent set of size
$W_1=(g-1)m/2$, and with $P=\{0,2,4\}$ (so $p=3$) gives a maximal independent set of size $W_3=(g-3)m/2$. The indices $0,2$ and $4$ are pairwise non-consecutive because $g\ge6$, and $p=3\le\lfloor g/2\rfloor$, so $P$ satisfies the hypotheses of Proposition~\ref{complete even}. Applying the proposition gives a maximal independent set of size $W_3=(g-3)m/2$. Then
$W_1-W_3=m>0$. In both cases two admissible cyclic sequences of different total sizes
exist, hence $G$  is not well-covered.
\end{proof}

\subsection{Non well-coveredness when $m$ is odd}
In this section, we study well-coveredness for odd $m$. Throughout this subsection $m\ge5$ is odd, and $(S_0,\dots,S_{g-1})$ is an admissible cyclic sequence, for shifts $\sigma_0,\dots,\sigma_{g-1}\in\mathbb{Z_m}$ and $g\ge3$. For odd $m$ the cycle $C_m$ has no pair of alternating maximum independent sets as in the even case, so we fix one maximum independent set $S_{max}=\{0,2,4,\dots,m-3\}, \ |S_{max}|=\tfrac{m-1}{2},$ and put $R=\mathbb{Z_m}\setminus S_{max}=\{1,3,\dots,m-2\}\cup\{m-1\}$, so that $|R|=\frac{m+1}{2}$. In $C_m[R]$ the vertices $1,3,\dots,m-2$ are pairwise non-adjacent and $m-1$ is adjacent to $m-2$ and to $0\notin R$; hence $C_m[R]$ has exactly one edge, $\{m-2,m-1\}$. Using $S_{max}$ or its translates in every block produces only sequences of one total size, so we also need an \emph{intermediate} independent set $T$ in some block, whose size can be varied.

\begin{lemma}\label{diff}
Let $S_{\max}=\{0,2,4,\ldots,m-3\}\subseteq\mathbb{Z}_m$, where $m\ge5$ is odd. Then
$\operatorname{Diff}(S_{\max})=\mathbb{Z}_m\setminus\{\pm1\}$, where $\operatorname{Diff}(S_{\max})=\{x-y:x,y\in S_{\max}\}$.
\end{lemma}

\begin{proof}
Write $S_{\max}=\{2i:0\le i\le r\}$, where $r=(m-3)/2$. Then every element of $\operatorname{Diff}(S_{\max})$ is of the form $2(i-j)$, where $-r\le i-j\le r$. Since $2r+1=m-2$, the set $\{-r,-r+1,\ldots,r\}$ contains exactly $m-2$ integers. Since $\gcd(2,m)=1$, multiplication by $2$ is a bijection on $\mathbb{Z}_m$;
therefore the values $2(i-j)$ produce exactly $m-2$ different elements of $\mathbb{Z}_m$. All possible values modulo $m$ can be represented by the integers $\{-(m-1)/2,\ldots,(m-1)/2\}$, and the set $\{-r,\ldots,r\}$ coincides with this representative system except for omitting the two values $\pm(m-1)/2$. To say explicitly, $\{-r,\ldots,r\}\cup\{\pm(m-1)/2\}=\{-(m-1)/2,\ldots,(m-1)/2\}$ and $\{-r,\ldots,r\}\cap\{\pm(m-1)/2\}=\emptyset$. Doubling these endpoints gives $2\cdot(m-1)/2=m-1\equiv-1\pmod m$ and $2\cdot(-(m-1)/2)=-(m-1)\equiv1\pmod m$, so $2\cdot\{\pm(m-1)/2\} =\{\pm1\}$. Since doubling is injective on the complete system $\{-(m-1)/2,\ldots,(m-1)/2\}$, it carries its complement $\{\pm(m-1)/2\}$ onto the complement of $\operatorname{Diff}(S_{\max})$ in $\mathbb{Z}_m$; hence the only residues missing from the difference set are $1$ and $-1$. Therefore $\operatorname{Diff}(S_{\max})=\mathbb{Z}_m\setminus\{\pm1\}$.
\end{proof}

\begin{lemma}
Let $S_{\max}+p$ and $S_{\max}+q$ be two translates of $S_{\max}$ where $p,q \in \mathbb{Z}_m$. Then $(S_{\max}+q)\cap(S_{\max}+p+\gamma)=\varnothing$ if and only if $q=p+\gamma+\varepsilon$, where $\varepsilon\in\{\pm1\}$.
\end{lemma}

\begin{proof}
The condition $(S_{\max}+q)\cap(S_{\max}+p+\gamma)=\varnothing$
holds if and only if $q-p-\gamma\notin\operatorname{Diff}(S_{\max})$. By the Lemma~\ref{diff},
$\operatorname{Diff}(S_{\max})=\mathbb{Z}_m\setminus\{\pm1\}$. Therefore $q-p-\gamma \in\{\pm1\}$.
Hence, there exists $\varepsilon\in\{\pm1\}$ such that $q-p-\gamma =\varepsilon$, or equivalently, $q=p+\gamma +\varepsilon$.

Conversely, if $q=p+\gamma +\varepsilon$, where $\varepsilon\in\{\pm1\}$, then $q-p-\gamma =\varepsilon\notin\operatorname{Diff}(S_{\max})$, again by the Lemma~\ref{diff}. Therefore $(S_{\max}+q)\cap(S_{\max}+p+\gamma )=\varnothing$.
\end{proof}

\begin{lemma}\label{DSP}
Let $S\subseteq\mathbb{Z}_m$ be an independent set. Then, for every $p\in\mathbb{Z}_m$, $D(S+p)=D(S)+p$. Moreover $D(S_{\max})=\emptyset$, hence $D(S_{\max}+p)=\emptyset$ for all $p$.
\end{lemma}

\begin{proof}
Recall that $D(S)=\mathbb{Z}_m\setminus N[S]$, where $N[S]$ denotes the closed neighborhood of $S$. Here $S+p=\{x+p\pmod m:x\in S\}.$ Therefore, $N[S+p]=N[S]+p$. Hence $D(S+p)=\mathbb{Z}_m\setminus N[S+p] =\mathbb{Z}_m\setminus(N[S]+p) =(\mathbb{Z}_m\setminus N[S])+p =D(S)+p$. Thus $D(S+p)=D(S)+p$.

The set $S_{\max}$ is a maximal independent set of $C_m$. Therefore, every vertex of $V(C_m)\setminus S_{\max}$ is adjacent to at least one vertex of $S_{\max}$. So every vertex lies in $N[S_{\max}]$ and $D(S_{\max})=\emptyset$.
\end{proof}

\begin{lemma} \label{T}
For $\delta\in\{0,1,\dots,m-1\}$ define $f(\delta)\in\mathbb{Z}$ by
\[
  f(\delta)=
  \begin{cases}
    \dfrac{m-1}{2}, & \delta=0,\\[2mm]
    \dfrac{\delta+1}{2}, & \delta\ \text{odd},\\[2mm]
    \dfrac{m-\delta+1}{2}, & \delta\neq0\ \text{even}.
  \end{cases}
\]
Let $u,v\in\mathbb{Z}_m$, put
$W_{u,v}=\mathbb{Z}_m\setminus\bigl((S_{\max}+u)\cup(S_{\max}+v)\bigr)$,
and let $T_{u,v}$ be a maximal independent set of the induced subgraph
$C_m[W_{u,v}]$. Then $T_{u,v}$ is an independent set of $C_m$,
$T_{u,v}\cap(S_{\max}+u)=\emptyset, T_{u,v}\cap(S_{\max}+v)=\emptyset,  D(T_{u,v})\subseteq(S_{\max}+u)\cup(S_{\max}+v),$ and $|T_{u,v}|=f(\delta), \text{where }\delta\in\{0,\dots,m-1\}\text{ is the remainder  of }
  u-v\pmod m .$
\end{lemma}

\begin{proof}
Since $T_{u,v}\subseteq W_{u,v}$ and $W_{u,v}$ is disjoint from both translates, $T_{u,v}$ misses both translates; and since $T_{u,v}$ is independent in the induced subgraph $C_m[W_{u,v}]$, it is independent in $C_m$. By maximality of $T_{u,v}$ in $C_m[W_{u,v}]$, every vertex of $W_{u,v}$ lies in the closed neighbourhood of $T_{u,v}$ in $C_m[W_{u,v}]$, so
$W_{u,v}\subseteq N_{C_m[W_{u,v}]}[T_{u,v}]=N_{C_m}[T_{u,v}]\cap W_{u,v}$, hence $W_{u,v}\subseteq N_{C_m}[T_{u,v}]$; taking complements in $\mathbb{Z}_m$ gives $D(T_{u,v})=\mathbb{Z}_m\setminus N_{C_m}[T_{u,v}]
  \subseteq\mathbb{Z}_m\setminus W_{u,v}
  =(S_{\max}+u)\cup(S_{\max}+v).$

Let $R=\mathbb{Z}_m\setminus S_{\max}=\{1,3,\dots,m-2\}\cup\{m-1\}$, so $|R|=(m+1)/2$. For the size, translate by $-v$. Then $|T_{u,v}|$ equals the size of a maximal independent set of $C_m[W']$, where $W'=W_{u,v}-v
    =\mathbb{Z}_m\setminus\bigl((S_{\max}+\delta)\cup S_{\max}\bigr)
    =R\cap(R+\delta)
    =\{x\in R:\ x-\delta\in R\},$
with $\delta\equiv u-v\pmod m$ and
$R=\mathbb{Z}_m\setminus S_{\max}=\{1,3,\dots,m-2\}\cup\{m-1\}$,
$|R|=(m+1)/2$. We first show that $C_m[W']$ has no edge when
$\delta\neq0$. Suppose $\{a,b\}\subseteq W'$ is an edge of $C_m$. Since $W'\subseteq R$ and $C_m[R]$ has exactly the edge $\{m-2,m-1\}$, we have $\{a,b\}=\{m-2,m-1\}$. Since $W'-\delta\subseteq R$ and translation preserves adjacency, $\{a-\delta,b-\delta\}$ is also an edge of $C_m[R]$, hence equals $\{m-2,m-1\}$. Comparing the two equalities yields $\{m-2-\delta,\,m-1-\delta\}=\{m-2,\,m-1\}$. If
$m-2-\delta=m-2$, then $\delta\equiv0$; if $m-2-\delta=m-1$ and
$m-1-\delta=m-2$, then $\delta\equiv-1$ and $\delta\equiv1$
simultaneously, so $2\equiv0\pmod m$, impossible because $m\ge5$ is odd. Hence $\delta\equiv0\pmod m$, contrary to the assumption $\delta\neq0$. Thus $C_m[W']$ has no edge for $\delta\neq0$, its only maximal independent set is $W'$ itself, and $|T_{u,v}|=|W'|$. We calculate the cardinality of  $W'=R\cap(R+\delta)$.

\emph{Case $\delta=0$.} Then $W'=R$, and every maximal independent set of $C_m[R]$ omits exactly one endpoint of $\{m-2,m-1\}$, so $|T_{u,v}|=|R|-1=\frac{m+1}{2}-1=\frac{m-1}{2}=f(0).$

\emph{Case $\delta\neq0$ even.} For odd $x\in\{1,3,\dots,m-2\}$ with $x>\delta$, the difference $x-\delta$ is odd and lies in
$\{1,3,\dots,m-2\}\subseteq R$, so such $x$ belong to $W'$; these are $x=\delta+1,\delta+3,\dots,m-2$, exactly $(m-\delta-1)/2$ elements. For odd $x<\delta$, the residue of $x-\delta$ modulo $m$ is $m+x-\delta$, which is even since $m$ and $x$ are odd and $\delta$ is even; the only even element of $R$ is $m-1$, so $m+x-\delta=m-1$, i.e.\ $x=\delta-1$: exactly one element. The remaining element $x=m-1$ gives $m-1-\delta$ even and different from $m-1$, hence not in $R$, so $x\notin W'$. Hence
$|T_{u,v}|=|W'|=\frac{m-\delta-1}{2}+1=\frac{m-\delta+1}{2}=f(\delta).$

\emph{Case $\delta$ odd.} For odd $x\ge\delta$ with $x\le m-2$, the difference $x-\delta$ is even and at most $m-3<m-1$, so it is not in $R$, and no such $x$ belongs to $W'$. For odd $x<\delta$, the residue $m+x-\delta$ is odd and lies in $\{1,3,\dots,m-2\}\subseteq R$; these are $x=1,3,\dots,\delta-2$, giving $(\delta-1)/2$ elements. The remaining element $x=m-1$ gives $m-1-\delta$ odd and lying in $\{1,3,\dots,m-2\}\subseteq R$, so $m-1\in W'$.  Hence $|T_{u,v}|=|W'|=\frac{\delta-1}{2}+1=\frac{\delta+1}{2}=f(\delta).$
In all three cases $|T_{u,v}|=f(\delta)$.
\end{proof}

\begin{lemma}\label{lem:preimages}
Let $m\ge5$ be odd, and let $f$ be as in Lemma~\ref{T}.
Then $f(\delta)=\dfrac{m-1}{2}$ if and only if
$\delta\in\{0,2,m-2\}$, and every value $v\neq\dfrac{m-1}{2}$ of $f$
has at most one odd and at most one even nonzero preimage.
Consequently, $f(\delta)<\dfrac{m-1}{2}$ for every
$\delta\notin\{0,2,m-2\}$; every value of $f$ has at most three
preimages; and four distinct residues
$\delta_1,\dots,\delta_4\in\{0,1,\dots,m-1\}$ cannot all have the same value of $f$.
\end{lemma}

\begin{proof}
  By Lemma~\ref{T}, we have $f(0)=\frac{m-1}{2},
  \ f(\delta)=\frac{\delta+1}{2}\ \ (\delta\text{ odd}),
  \ f(\delta)=\frac{m-\delta+1}{2}\ \ (\delta\text{ even},\ \delta\neq0).$ These three formulas cover the three disjoint classes $\{0\}$, the odd residues, and the even nonzero residues. Note that $m-2$ is odd and nonzero and that $2$ is even and nonzero, since $m$ is odd and $m\ge5$.

 Calculating the preimage of others values apart from $\frac{m-1}{2}$ using the values of function $f(\delta) $. The residue $0$ is a preimage by the first formula. On the odd residues the function $\delta\mapsto\frac{\delta+1}{2}$ is strictly increasing, hence injective, and it equals $\frac{m-1}{2}$ if and only if $\delta=m-2$. On the even nonzero residues the function $\delta\mapsto\frac{m-\delta+1}{2}$ is strictly decreasing, hence injective, and it equals $\frac{m-1}{2}$ if and only if $\delta=2$. Since the three classes are disjoint and each contains exactly one preimage, we obtain
$f^{-1}\!\left(\frac{m-1}{2}\right)=\{0,2,m-2\}$, exactly three
preimages. Now consider the preimages apart from $\frac{m-1}{2}$, consider $v\neq\frac{m-1}{2}$. Since $f(0)=\frac{m-1}{2}$, the residue $0$ is not a preimage of $v$. On the odd residues $f$ is injective, so $v$ has at most one odd preimage; on the even nonzero residues $f$ is injective, so $v$ has at most one even nonzero preimage. Hence $v$ has at most two preimages in total.

Next to show $f(\delta)<\frac{m-1}{2}$, by Lemma~\ref{T} we have
$f(0)=\frac{m-1}{2}$,
$f(\delta)=\frac{\delta+1}{2}\le\frac{m-1}{2}$ for odd
$\delta\le m-2$, with equality if and only if $\delta=m-2$, and
$f(\delta)=\frac{m-\delta+1}{2}\le\frac{m-1}{2}$ for even $\delta\ge2$, with equality if and only if $\delta=2$. Hence
$f(\delta)\le\frac{m-1}{2}$ for every $\delta$, with equality exactly for $\delta\in\{0,2,m-2\}$; therefore
$f(\delta)<\frac{m-1}{2}$ for every $\delta\notin\{0,2,m-2\}$.

 The value $\frac{m-1}{2}$ has exactly three
preimages, and every other value has at most two; hence every value of $f$ has at most three preimages. In particular, four distinct residues cannot share the same value of $f$.
\end{proof}

\begin{lemma}\label{lem:one-intermediate}
Let $m\ge 5$ be odd and $g\ge 3$, let $\sigma\in\mathbb{Z}_m^g$ and $\Lambda=\sum\limits_{j\in\mathbb{Z}_g}\sigma_j$. Fix a block $k\in\mathbb{Z}_g$.
For each $j\in\mathbb{Z}_g\setminus\{k-1,k\}$, choose $\epsilon_j\in\{\pm 1\}$ and put $E'=\sum\limits_{j\notin\{k-1,k\}}\epsilon_j$. Then there exists an admissible cyclic sequence $(S_0,\dots,S_{g-1})$ such that $S_j=S_{\max}+p_j$ for every $j\neq k$ and $|S_k|=f(\delta)$, where $\delta\equiv\Lambda+E'\pmod m$.
\end{lemma}

\begin{proof}
We start by fixing the phase of the block immediately after $k$ by setting $p_{k+1}=0$. We now walk around the path $k+1 \to k+2 \to \cdots \to k-1$ . For each step $j$ along this path, for the translates of $S_{\max}$  to avoid edges between block $j$ and block $j+1$, their phases must satisfy $p_{j+1} \equiv p_j + \sigma_j \pm 1 \pmod m$. The freedom to choose $+1$ or $-1$ is exactly what the terms $\epsilon_j$ represent. Thus, we define the phases recursively by $p_{j+1} = p_j + \sigma_j + \epsilon_j.$ Because every translate of $S_{\max}$ has an empty defect set, all maximality conditions along this path are automatically satisfied. Summing these relations from $j=k+1$ to $k-2$ gives the total phase difference $p_{k-1}-p_{k+1} \equiv \sum\limits_{j\notin\{k-1,k\}}\sigma_j + E' \pmod m.$ Since $p_{k+1}=0$, this determines the phase $p_{k-1}$. To place a state in block $k$, we need to know what it must avoid. From the perspective of block $k$, the left neighbour $S_{k-1}$, when shifted by $\sigma_{k-1}$, looks like $S_{\max} + u$ where $u = p_{k-1} + \sigma_{k-1}$. Similarly, the right neighbour $S_{k+1}$, when pulled back by $\sigma_k$, looks like $S_{\max} + v$ where $v = p_{k+1} - \sigma_k$. By Lemma~\ref{T}, applied to these $u$ and $v$, there exists a set $T_{u,v}$ that avoids both $S_{\max}+u$ and $S_{\max}+v$, and whose undominated vertices are contained in their union. We set $S_k = T_{u,v}$.

The two disjointness conditions of $T_{u,v}$ guarantee independence across the transitions $(k-1, k)$ and $(k, k+1)$. The containment of the defect set guarantees maximality at block $k$. All other blocks carry translates of $S_{\max}$, so their defect sets are empty and maximality holds automatically. Thus, $(S_0, \dots, S_{k-1}, T_{u,v}, S_{k+1}, \dots, S_{g-1})$ is admissible. Finally, the parameter $\delta$ controlling the size of $T_{u,v}$ is given by, $\delta \equiv u-v = p_{k-1}-p_{k+1}+\sigma_{k-1}+\sigma_k \equiv \sum_{j\notin\{k-1,k\}}\sigma_j + E' + \sigma_{k-1} + \sigma_k \equiv \Lambda+E' \pmod m.$ Thus, $|S_k| = |T_{u,v}| = f(\delta)$.
\end{proof}

\begin{remark}\label{rem:intermediate-size}
The formula $\delta\equiv\Lambda+E'\pmod m$ shows that the intermediate block attains the maximal possible size $f(\delta)=\frac{m-1}{2}$ exactly when $\delta\in\{0,2,m-2\}$.  To understand why, recall that $f(\delta)$ is maximized exactly when the intermediate block is as large as a translate of $S_{\max}$. If we had placed a translate of $S_{\max}$ in block $k$ as well, the two transitions incident to $k$ would have corrections $\epsilon_{k-1},\epsilon_k\in\{\pm 1\}$. Their sum can only be $-2, 0$, or $2$. The cycle would close perfectly as an all-maximum-state sequence if and only if $\Lambda+E'+\epsilon_{k-1}+\epsilon_k \equiv 0 \pmod m,$ which is solvable exactly when $\Lambda+E'\in\{0,2,m-2\}\pmod m$, i.e., when $\delta\in\{0,2,m-2\}$. Thus, the intermediate block fails to be maximal (i.e., $|S_k|<\frac{m-1}{2}$) precisely when the two missing corrections $\epsilon_{k-1}$ and $\epsilon_k$ \emph{cannot} close the cycle. This strict size drop is exactly what we need to prove that the graph is not well-covered.
\end{remark}

\begin{proposition}\label{prop:g5}
Let $m\ge5$ be odd and $g\ge5$. Then $G=C_n(a,b)$ is not
well-covered.
\end{proposition}

\begin{proof}
To prove that $G$ is not well-covered, it suffices to exhibit two maximal independent sets of different sizes. By Lemma~\ref{lem:one-intermediate}, we can build an independent set by choosing corrections $\epsilon_j \in \{\pm 1\}$ on all transitions except the two incident to a fixed block $k$. The total size of the resulting independent set is $(g-1)\frac{m-1}{2} + f(\delta)$, where $\delta\equiv\Lambda+E'\pmod m$ and $E'=\sum \epsilon_j$.

Since $g\ge5$, the number of available transitions not incident to block $k$ is $r=g-2\ge3$. This gives us enough freedom to construct four different choices for the corrections. For each $t\in\{0,1,2,3\}$, we choose exactly $t$ of the $r$ corrections to be $+1$ and the remaining $r-t$ to be $-1$. Thus, the total corrections are calculated as  $E'_t = t(+1) + (r-t)(-1) = 2t-r.$ The four resulting values are $E'_0=-r$, $E'_1=-r+2$, $E'_2=-r+4$, and $E'_3=-r+6$. The pairwise differences between these values are $2, 4$, and $6$. Since $m\ge5$ is odd, $m$ divides none of these differences, meaning the four values $E'_t$  are pairwise distinct modulo $m$ and hence the four values of $\delta_t\equiv\Lambda+E'_t \pmod m$. For each $t$, Lemma~\ref{lem:one-intermediate} provides an admissible cyclic sequence of size  $W_t=(g-1)\frac{m-1}{2}+f(\delta_t).$ By Lemma~\ref{lem:preimages}, the function $f$ cannot take the same value at four distinct residues. Thus, among our four distinct values $\delta_0,\delta_1,\delta_2,\delta_3$, the values $f(\delta_t)$ cannot all be equal. Consequently, at least two of the sequences will have different total sizes $W_t \neq W_{t'}$. Hence, $G$ is not well-covered.
\end{proof}

\begin{proposition}\label{prop:g4}
Let $m\ge5$ be odd and $g=4$. Then $G=C_n(a,b)$ is not well-covered.
\end{proposition}

\begin{proof}
To prove that $G$ is not well-covered, we again look for two maximal independent sets of different sizes. Here, $g=4$, so the four-shot method from Proposition~\ref{prop:g5} fails here as there are only $g-2=2$ available transitions, giving at most three distinct correction values. Instead, we will use two intermediate blocks and leave the other two blocks as translates of $S_{\max}$.

Put $\Lambda=\sigma_0+\sigma_1+\sigma_2+\sigma_3\pmod m$. Choose $p_0,p_2\in\mathbb{Z_m}$ and set $S_0=S_{\max}+p_0$, $S_2=S_{\max}+p_2$. Apply Lemma~\ref{T} twice, first with $u_1=p_0+\sigma_0$ and $v_1=p_2-\sigma_1$ to obtain an intermediate state $T_1$, and then with $u_2=p_2+\sigma_2$ and $v_2=p_0-\sigma_3$ to obtain an intermediate state $T_3$. Then $(S_0,T_1,S_2,T_3)$ is admissible, the independence conditions at the four transitions and the maximality condition at $C_1$ and $C_3$ follow from Lemma~\ref{T} exactly as in Lemma~\ref{lem:one-intermediate}, while maximality at $C_0$ and
$C_2$ holds because $D(S_{\max}+p)=\emptyset$ for every $p\in\mathbb{Z_m}$. The parameters of the two intermediate states are $\delta_1\equiv u_1-v_1=p_0-p_2+\sigma_0+\sigma_1\pmod m, \delta_2\equiv u_2-v_2=p_2-p_0+\sigma_2+\sigma_3\pmod m,$ so that $\delta_1+\delta_2\equiv\Lambda\pmod m$. The total size of
the corresponding maximal independent set is $  W=(m-1)+f(\delta_1)+f(\Lambda-\delta_1)=(m-1)+h(\delta_1),$ where $h(\delta)=f(\delta)+f(\Lambda-\delta)$ for $\delta\in\{0,1,\dots,m-1\}$. Since $p_0,p_2$ are arbitrary,
$\delta_1$ ranges over all of $\mathbb{Z_m}$, so it suffices to show that $h$ is not constant.

Let $\Lambda\in\{0,1,\dots,m-1\}$. We distinguish five cases, using the formula for $f$ in Lemma~\ref{T}.
\emph{Case $\Lambda=0$.} Then
$h(0)=f(0)+f(0)=m-1$ and $h(1)=f(1)+f(m-1)=1+1=2$; since $m\ge5$,
$m-1\neq2$.

\emph{Case $\Lambda$ even, $2\le\Lambda\le m-3$.} Then $ h(0)=f(0)+f(\Lambda)=\frac{m-1}{2}+\frac{m-\Lambda+1}{2}
       =m-\frac{\Lambda}{2}, h(1)=f(1)+f(\Lambda-1)=1+\frac{\Lambda}{2}.$ These coincide only when $\Lambda=m-1$, which the assumed range excludes.

\emph{Case $\Lambda$ odd, $3\le\Lambda\le m-2$.} Then $h(0)=f(0)+f(\Lambda)=\frac{m-1}{2}+\frac{\Lambda+1}{2}
       =\frac{m+\Lambda}{2}, h(1)=f(1)+f(\Lambda-1)=1+\frac{m-\Lambda+2}{2} =\frac{m-\Lambda+4}{2},$ whose difference $h(0)-h(1)=\Lambda-2$ is nonzero because $\Lambda\ge3$.

\emph{Case $\Lambda=1$.} Then $h(0)=f(0)+f(1)=\frac{m-1}{2}+1$ and $h(3)=f(3)+f(m-2)=2+\frac{m-1}{2}$, so $h(3)-h(0)=1\neq0$.

\emph{Case $\Lambda=m-1$.} Then $h(0)=f(0)+f(m-1)=\frac{m-1}{2}+1$ and $h(m-3)=f(m-3)+f(2)=2+\frac{m-1}{2}$, so $h(m-3)-h(0)=1\neq0$.

These five cases cover every $\Lambda\in\{0,1,\dots,m-1\}$, and in each case we find two values $\delta,\delta'$ with
$h(\delta)\neq h(\delta')$. The corresponding admissible sequences have total sizes $W(\delta)=(m-1)+h(\delta)$ and
$W(\delta')=(m-1)+h(\delta')$, which are different. So, $G$ is not well-covered.
\end{proof}

\begin{theorem}\label{wellcovered}
Let $n,a,b\in\mathbb{Z}_{>0}$ with $1\leq a<b\leq\left\lfloor\frac{n-1}{2}\right\rfloor$ and $\gcd(n,a,b)=1$ with $n\geq19$, and let $ G =C_n(a,b)$.
\begin{enumerate}
     \item[(a)]  $G$ is well-covered with $\alpha( $G $)=n/3$ if $n=3a$ or $n=3b$ for $n>6$.
    \item[(b)]  If $n\neq3a$ and $n\neq3b$, then $G$ is not well-covered, provided   $ \begin{cases}   m_c \geq4 \text{    even}, & g_c \geq6,\\[4pt]   m_c\geq5  \text{    odd}, & g_c\geq4\\[4pt]   \end{cases} $.
\end{enumerate}

\end{theorem}
\begin{proof}
\begin{enumerate}
    \item[(a)] Follows from Theorem~ \ref{thm:well-covered}.
    \item[(b)] If $n\neq3b$, provided $m_c ( \geq4 ) \text{    even},  g_c \geq6$ then $G$ is not well-covered follows by the  Theorem~\ref{even case}.Next, if $n\neq3b$, provided $m_c ( \geq5 ) \text{ odd},  g_c \geq4$, then $G$ is not well-covered follows from the Propositions \ref{prop:g5}\ref{prop:g4}. For  $n\neq3 a$ the result follows similarly.
\end{enumerate}
    \end{proof}

\section{\centering NON-COHEN MACAULAYNESS OF $C_n(a,b)$}
In this section, we prove that $C_n(a,b)$ is not Cohen- Macaulay when $n=3b$ and $3a$. Here we  construct a suitable face $F$ of the independence complex in such a manner that the link $\operatorname{lk}_{\Delta}(F)$ is disconnected. Then, we apply Reisner's criterion.

\begin{lemma}\label{lem:link-edge-criterion}
Let $F\in\Delta$, and let $u,v\in V(\operatorname{lk}_{\Delta}(F))$. Then $\{u,v\}\in\operatorname{lk}_{\Delta}(F)$ if and only if $u$ and $v$ are non-adjacent in $G$.
\end{lemma}

\begin{proof}
Since $u,v\in V(\operatorname{lk}_{\Delta}(F))$, by the definition of the link, both $F\cup\{u\}$ and $F\cup\{v\}$ are independent sets in $G$. Thus, neither $u$ nor $v$ is adjacent to any vertex of $F$.Therefore, when the two vertices $u$ and $v$ are considered together, the only possible new adjacency that can prevent $F\cup\{u,v\}$ from being independent is an adjacency between $u$ and $v$. Hence,
$F\cup\{u,v\}$ is independent if and only if $u$ and $v$ are non-adjacent in $G$. Also the link, $\{u,v\}\in\operatorname{lk}_{\Delta}(F)$ if and only if
$F\cup\{u,v\}\in\Delta$. Since $\Delta=\operatorname{Ind}(G)$, this is equivalent to $F\cup\{u,v\}$ being an independent set in $G$. Therefore, $\{u,v\}\in\operatorname{lk}_{\Delta}(F)$ if and only if $u$ and $v$
are non-adjacent in $G$.
\end{proof}

\begin{lemma} \label{reduced-zero-homology}\cite{miller2004combinatorial}
Let $X$ be a nonempty simplicial complex with $q$ connected components. Then, for every field $K$, $\dim_K\widetilde{H}_0(X;K)=q-1$.
\end{lemma}
\vspace{0.2cm}
\begin{theorem}\cite{reisner1976cohen}[Reisner's Criterion] \label{reisner}
Let $\Delta$ be a simplicial complex and let $K$ be a field. Then $\Delta$ is Cohen--Macaulay over $K$ if and only if $\widetilde{H}_i(\operatorname{lk}_{\Delta}(F);K)=0$ for every face $F\in\Delta$ and every integer $i$ satisfying $i<\dim(\operatorname{lk}_{\Delta}(F))$.

\end{theorem}

\begin{theorem}\cite{Milnor1956}[K\"unneth formula]\label{kunneth}
Let $X$ and $Y$ be simplicial complexes on disjoint vertex sets and let $\mathbb{K}$ be a field. Then for every integer $j$, $\widetilde{H}_j(X*Y;\mathbb{K})\;\cong\;\bigoplus_{p+q=j-1}\widetilde{H}_p(X;\mathbb{K})\otimes_{\mathbb{K}}\widetilde{H}_q(Y;\mathbb{K}),$ where $\widetilde{H}_*$ denotes reduced simplicial homology, with the convention $\widetilde{H}_{-1}(\{\emptyset\};\mathbb{K})=\mathbb{K}$.
\end{theorem}
\subsection{Setup and Structural Description for $n=3b$}

Let $n,a,b\in\mathbb{Z}_{>0}$, with $1\leq a<b\leq \left\lfloor\frac{n-1}{2}\right\rfloor$, and suppose that
$n=3b$, where $b\geq 5$ and $\gcd(n,a,b)=1$. Let $G=C_n(a,b)$ and let $\Delta=\operatorname{Ind}(G)$ denote the independence complex of $G$.

We now describe the structure of $G$ that is used in the construction of the face $F$. Since $n=3b$, every vertex $v\in\mathbb{Z}_n$ has a unique representation
$v=t+cb$, where $0\leq t\leq b-1$ and $c\in\{0,1,2\}$. We denote the vertex $t+cb$ by $T_t^c$. For each $t\in\{0,\ldots,b-1\}$, let $T_t=\{T_t^0,T_t^1,T_t^2\}$. Thus, $T_t^0=t$, $T_t^1=t+b$, and $T_t^2=t+2b$, and the vertex set of $G$ is partitioned as $V(G)=T_0\cup T_1\cup\cdots\cup T_{b-1}$. The edges corresponding to the generator $b$ are considered, for every $t$,  $T_t^1-T_t^0=b$ and $T_t^2-T_t^1=b$, while $T_t^0-T_t^2=-2b\equiv b\pmod{3b}$. Since the connection set contains $\pm b$, the three vertices in
$T_t$ are pairwise adjacent. Therefore, $G[T_t]\cong K_3$ for every $t$. Hence the vertices of $G$ are partitioned into the $b$ triangles $T_0,T_1,\ldots,T_{b-1}$. In particular, any independent set of $G$ can contain at most one vertex from each triangle $T_t$.

It remains to understand how the edges corresponding to $a$ connect these triangles. For a vertex $T_t^c=t+cb$, addition by $a$ gives $T_t^c+a=t+a+cb$. The resulting expression depends on whether
$t+a<b$ or $t+a\geq b$. If $t+a<b$, then $T_t^c+a=T_{t+a}^c$. On the other hand, if $t+a\geq b$, then
$t+a=(t+a-b)+b$, and hence $T_t^c+a=(t+a-b)+(c+1)b=T_{t+a-b}^{c+1}$, where the second superscript is taken modulo $3$. For convenience, define $\varepsilon_t=
\begin{cases}
0,&t+a<b,\\
1,&t+a\geq b.
\end{cases}$

Then the two cases can be written uniformly as \begin{equation}\label{eq:n3b-edge-action}
T_t^c+a
=
T_{(t+a)\bmod b}^{\,c+\varepsilon_t\pmod 3}.
\end{equation}

Thus, equation~\eqref{eq:n3b-edge-action} describes precisely how an edge corresponding to $a$ connects a vertex in one triangle to a vertex in another triangle.

Since $\gcd(n,a,b)=1$, repeatedly adding $a$ modulo $b$ visits every element of $\mathbb{Z}_b$ exactly once. Define
$q_k=ka\pmod b$, $0\leq k\leq b-1$. Then $q_0,q_1,\ldots,q_{b-1}$ are all the elements of $\mathbb{Z}_b$,
and $q_{k+1}\equiv q_k+a\pmod b$. Therefore, the triangles can be arranged in cyclic order
$T_{q_0},T_{q_1},\ldots,T_{q_{b-1}},T_{q_0}$. More precisely, for each $k$, the map ${T_q^k}_c \mapsto {T_q^{k+1}}_{c + \varepsilon_q^k}$ is a bijection $T_q^k \to T_q^{k+1}$; hence the $a$-edges between consecutive triangles form a perfect matching, and every vertex of $T_q^k$ has exactly one $a$-neighbor in $T_q^{k+1}$ and exactly one in $T_q^{k-1}$. The $a$-edges join vertices only between two consecutive triangles in this order. Thus, the $a$-edges give a cycle through all the $b$ triangles.

\begin{figure}[!ht]
\centering
\begin{tikzpicture}[
  v/.style={circle,draw,fill=white,inner sep=0pt,minimum size=6.5mm,font=\scriptsize},
  trilabel/.style={font=\small,gray}]

  \def\RadO{4.4}\def\RadM{3.3}\def\RadI{2.2}

  \foreach \kk in {0,...,5}{
    \pgfmathsetmacro{\thk}{90-60*\kk}
    \pgfmathsetmacro{\thm}{90-60*\kk+14}
    \node[v] (O\kk) at (\thk:\RadO) {$T^{0}_{q_{\kk}}$};
    \node[v] (M\kk) at (\thm:\RadM) {$T^{1}_{q_{\kk}}$};
    \node[v] (I\kk) at (\thk:\RadI) {$T^{2}_{q_{\kk}}$};
    \draw[thick] (O\kk)--(M\kk)--(I\kk);
    \draw[thick] (I\kk) to[bend right=45] (O\kk);
    \node[trilabel] at (\thk:5.2) {$T_{q_{\kk}}$};}

  \foreach \kk/\nn in {0/1,1/2,2/3,3/4,4/5}{
    \foreach \rg in {O,M,I}{
      \draw[blue,dashed,thick] (\rg\kk)--(\rg\nn);
    }
  }

  \draw[red,dashed,thick] (O5) to[bend left=12] (M0);
  \draw[red,dashed,thick] (M5) to[bend left=12] (I0);
  \draw[red,dashed,thick] (I5) to[bend left=12] (O0);
\end{tikzpicture}
\caption{Block decomposition of $C_{3b}(a,b)$, drawn for $b=6$ and $a=1$ (so $n=18$ and $q_k=k$). Each $T_{q_k}$ is a triangle formed by the $b$-edges (solid black), with layers $T^{0}_{q_k},T^{1}_{q_k},T^{2}_{q_k}$. The $a$-edges join consecutive triangles $T_{q_k}$ and $T_{q_{k+1}}$ by a perfect matching, so together they form a cycle through all $b$ triangles. Between $T_{q_k}$ and $T_{q_{k+1}}$ ($0\le k\le4$) the matching preserves the layer ($\varepsilon_{q_k}=0$, blue dashed lines). The matching between the last triangle $T_{q_5}$ and the first one $T_{q_0}$ shifts the layer by $+1$ ($\varepsilon_{q_5}=1$, because $q_5+a=6\ge b$); it is drawn in red dashed lines and closes the cycle.}
\end{figure}
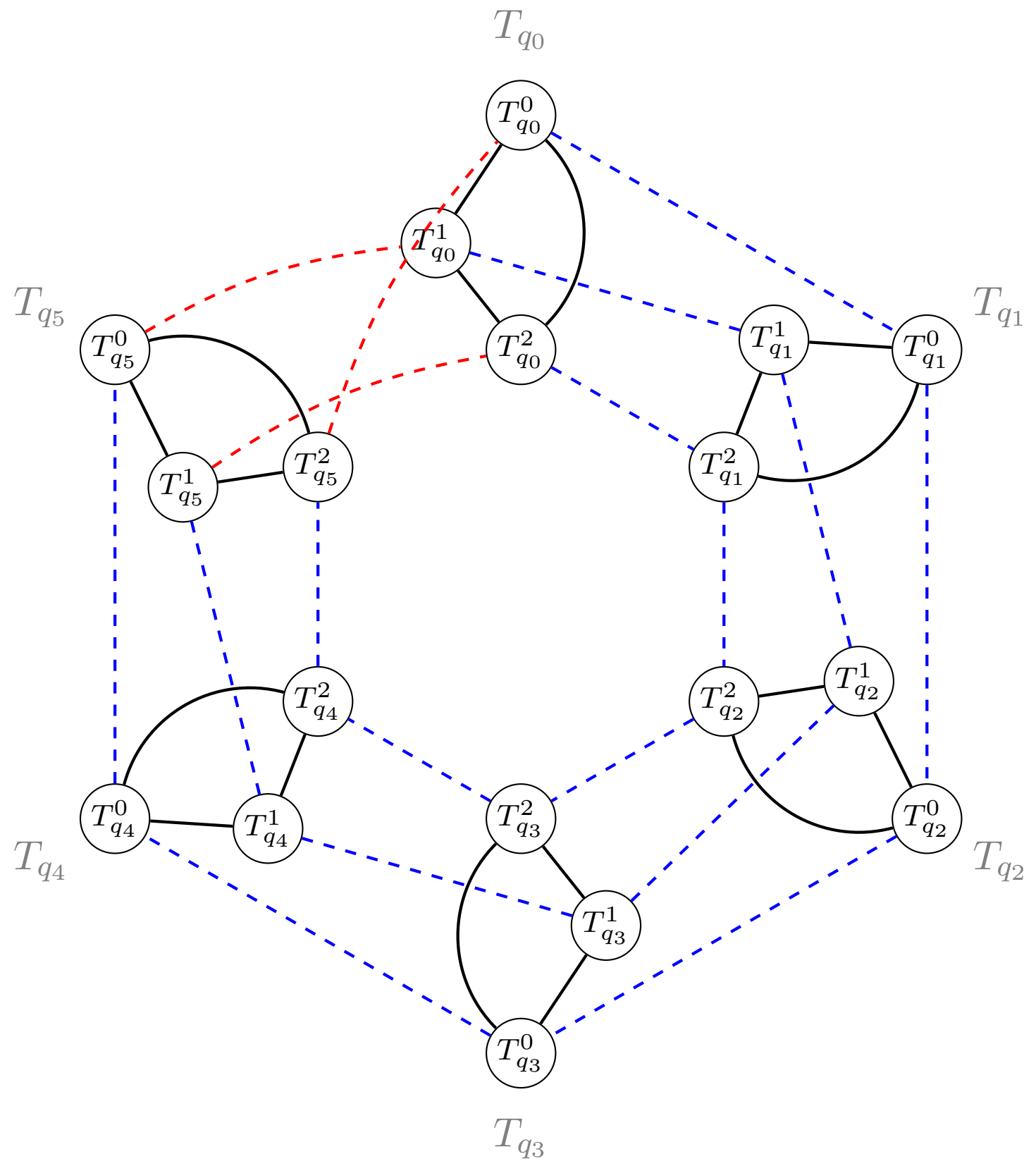

\subsection{Construction of the Face $F$}

We construct $F$ with two purposes. First, we want to choose one vertex from each triangle except two consecutive triangles, in such a way that the chosen vertices form an independent set.

Since $\gcd(n,a,b)=1$, the numbers $q_k\equiv ka\pmod b$, $0\leq k\leq b-1$, are distinct and give all the elements of $\mathbb{Z}_b$. Thus, the triangles occur in the cyclic order $T_{q_0},T_{q_1},\ldots,T_{q_{b-1}},T_{q_0}$, and consecutive triangles in this order are joined by $a$-edges. We omit the two consecutive triangles $T_{q_0}=T_0$ and $T_{q_1}=T_a$. From every remaining triangle $T_{q_k}$, $2\leq k\leq b-1$, we choose one vertex. Write the chosen vertex as $T_{q_k}^{c_k}$, where $c_k\in\{0,1,2\}$ here $q_k$ is represented as first coordinate and $c_k$ is represented as second coordinate. We now explain how the numbers $c_k$ are chosen. Suppose that $T_{q_k}^{c_k}$ has been chosen. From the formula for the $a$-edges, $T_{q_k}^{c_k}+a = T_{q_{k+1}}^{c_k+\varepsilon_{q_k}\pmod 3}.$ Thus, among the three vertices of $T_{q_{k+1}}$, exactly one is adjacent to $T_{q_k}^{c_k}$. Its second coordinate is
$c_k+\varepsilon_{q_k}$. Therefore, to prevent an $a$-edge between the two chosen vertices, we require $c_{k+1}\not\equiv c_k+\varepsilon_{q_k}\pmod 3.$ Next,  We want the same second coordinate to be excluded from the two omitted triangles. This property will be used later to determine the link.  For every $c_k$, there are two choices for $c_{k+1}$ satisfying this condition. We use this freedom to obtain the above property. Set $c_2=0$. For $2\leq k\leq b-2$, write $c_{k+1}\equiv c_k+\varepsilon_{q_k}+\delta_k\pmod 3,$ where $\delta_k\in\{1,2\}$. Since $\delta_k\neq0$, which gives $c_{k+1}\not\equiv c_k+\varepsilon_{q_k}\pmod 3$. We next choose the numbers $\delta_k$ so that the same second coordinate is excluded from $T_0$ and $T_a$. The last chosen vertex $T_{q_{b-1}}^{c_{b-1}}$ excludes from $T_0$ the vertex whose second coordinate is
$c_{b-1}+\varepsilon_{q_{b-1}}$. Similarly, the first chosen vertex $T_{q_2}^{c_2}$ excludes from $T_a$
the vertex whose second coordinate is $c_2-\varepsilon_{q_1}$.Therefore, we want $c_{b-1}+\varepsilon_{q_{b-1}}
\equiv c_2-\varepsilon_{q_1}\pmod 3.$ Put $\Lambda=\sum\limits_{k=2}^{b-2}\varepsilon_{q_k}.$ There are $b-3$ numbers $\delta_k$. Suppose exactly $U$ of them are equal to $2$, and the remaining $b-3-U$ are equal to $1$. Then $\sum\limits_{k=2}^{b-2}\delta_k=(b-3)+U.$ Using the above conditions,we have $\Lambda+(b-3)+U+\varepsilon_{q_{b-1}}
+\varepsilon_{q_1}\equiv0\pmod 3.$

Hence we choose $U\in\{0,1,2\}$ such that
\begin{equation}
U\equiv
-\left(
\Lambda+(b-3)+\varepsilon_{q_{b-1}}
+\varepsilon_{q_1}
\right)
\pmod 3.
\tag{1}
\end{equation}

Since $b\geq5$, we have $b-3\geq2$. Finally, define
\begin{equation}
F=\{T_{q_k}^{c_k}:2\leq k\leq b-1\}.
\tag{2}
\end{equation}

Thus, $F$ contains exactly one vertex from each of $T_{q_2},T_{q_3},\ldots,T_{q_{b-1}}$, and no vertex from $T_0$ or $T_a$. The choice of the numbers $c_k$ ensures that no two selected vertices in
consecutive triangles $T_{q_k}$ and $T_{q_{k+1}}$ ($2\le k\le b-2$) are
joined by an $a$-edge. The pair $(T_{q_{b-1}}, T_{q_2})$ is not consecutive
in the cyclic order, because $b\ge5$ implies $b-1\ge4$, so the two triangles
are separated by $T_{q_0}$ and $T_{q_1}$ on one side and by
$T_{q_3},\dots,T_{q_{b-2}}$ on the other; hence no $a$-edge joins a vertex
of $T_{q_{b-1}}$ to a vertex of $T_{q_2}$. Since at most one vertex is
chosen from each triangle, and $b$-edges stay inside triangles, no two
selected vertices are joined by a $b$-edge. Hence $F$ is independent, so
$F\in\Delta$, and $|F|=b-2$.
\begin{lemma}\label{surviving vertices}
Let $F$ be the face constructed above, and define $X\equiv c_{b-1}+\varepsilon_{q_{b-1}}\pmod{3}$ and
$Y\equiv c_2-\varepsilon_{q_1}\pmod{3}$. Then $X=Y$, and consequently $V(\operatorname{lk}_{\Delta}(F))
= \{T_0^i:i\neq X\}\cup\{T_a^i:i\neq X\}$.
\end{lemma}

\begin{proof}

Every triangle
$T_{q_k}$, $2\le k\le b-1$, contains the vertex $T_{q_k}^{c_k}\in F$;
since $T_{q_k}$ is a clique, every other vertex of $T_{q_k}$ is adjacent
to $T_{q_k}^{c_k}$ and therefore cannot belong to the link of $F$. The
vertices of $F$ themselves are not in the link. Hence $V(\operatorname{lk}_\Delta(F))\subseteq T_{q_0}\cup T_{q_1}=T_0\cup T_a.$ Our purpose is therefore to determine which vertices of these two triangles are excluded by the vertices of $F$. The last vertex chosen in $F$ is $T_{q_{b-1}}^{c_{b-1}}$. Its $+a$-neighbor lies in $T_0$ and has second coordinate
$X\equiv c_{b-1}+\varepsilon_{q_{b-1}}\pmod{3}$. Thus, $T_0^X$ is adjacent to a vertex of $F$ and cannot belong to the link. Similarly, the first vertex chosen in $F$ is $T_{q_2}^{c_2}$. Its
$-a$-neighbor lies in $T_a$ and has second coordinate $Y\equiv c_2-\varepsilon_{q_1}\pmod{3}$.
Thus, $T_a^Y$ is adjacent to a vertex of $F$ and cannot belong to the link. We now show that $X=Y$. From the recurrence $c_{k+1}\equiv c_k+\varepsilon_{q_k}+\delta_k\pmod{3}$,
summing for $2\leq k\leq b-2$ gives $c_{b-1}-c_2 \equiv \sum\limits_{k=2}^{b-2}\varepsilon_{q_k}
+
\sum\limits_{k=2}^{b-2}\delta_k
\pmod{3}$. By the choice of the $\delta_k$'s, $\sum_{k=2}^{b-2}\delta_k=(b-3)+U$. Hence $c_{b-1}-c_2\equiv
\Lambda+(b-3)+U \pmod{3}$. Therefore, $X-Y \equiv c_{b-1}-c_2 +\varepsilon_{q_{b-1}} +\varepsilon_{q_1} \pmod{3} \equiv
\Lambda+(b-3)+U +\varepsilon_{q_{b-1}} +\varepsilon_{q_1} \pmod{3} \equiv 0\pmod{3},$ by the definition of $U$. Since $X,Y\in\{0,1,2\}$, it follows that
$X=Y$. Thus, the same second coordinate is excluded from both $T_0$ and
$T_a$. Since each of these triangles contains three vertices, exactly
two vertices survive in each. Hence
$V(\operatorname{lk}_{\Delta}(F))
=
\{T_0^i:i\neq X\}\cup\{T_a^i:i\neq X\}$.
\end{proof}
\begin{lemma}\label{Structure of the Link}
Let $F$ be the face constructed above, and let $S=V(\operatorname{lk}_{\Delta}(F))$. Then $\operatorname{lk}_{\Delta}(F)\cong K_2\sqcup K_2$.
\end{lemma}

\begin{proof}
By Lemma~\ref{surviving vertices}, $S=\{T_0^i:i\neq X\}\cup\{T_a^i:i\neq X\}$. Thus the same second coordinate $X$ is excluded from both $T_0$ and $T_a$, and exactly two vertices survive in each triangle.

We now determine all adjacencies between $T_0$ and $T_a$.

For $i,j\in\{0,1,2\}$, $T_0^i=ib$ and $T_a^j=a+jb$. An edge corresponding to $+a$ gives
$T_0^i+a=ib+a=T_a^i$. Hence a $+a$-edge joins $T_0^i$ and $T_a^j$ precisely when $i=j$. Next, we consider a possible $-a$-edge. Suppose that
$T_0^i$ and $T_a^j$ are joined by a $-a$-edge. Then $T_a^j-T_0^i=a+(j-i)b\equiv -a\pmod{3b},$
and hence $2a\equiv(i-j)b\pmod{3b}.$ Set $k=i-j$. Since $i,j\in\{0,1,2\}$, we have
$k\in\{-2,-1,0,1,2\}$. Thus $2a-kb=3bq$ for some $q\in\mathbb{Z}$. Since $0<a<b$, we have
$0<2a<2b$. We consider the possible values of $k$. If $k=0$, then $2a=3bq$, which is impossible because $0<2a<3b$. If $k=1$, then
$-b<2a-b<b$, so the only multiple of $3b$ that $2a-b$ can equal is $0$, giving $2a=b$. If $k=-1$, then $b<2a+b<3b$, so $2a+b$ cannot be a multiple of $3b$. If $k=2$, then
$-2b<2a-2b<0$, so again no multiple of $3b$ is possible. Finally, if $k=-2$, then
$2b<2a+2b<4b$, so the only possible multiple of $3b$ is $3b$, which again gives $2a=b$. Thus, in every possible case, the congruence implies $2a=b$. Hence $b\mid2a$. Since $\gcd(n,a,b)=1$, we obtain $b\mid2$, contradicting
$b\ge5$. Therefore, no $-a$-edge joins $T_0$ and $T_a$. Finally, the $b$-edges lie entirely within $T_0$ or within $T_a$,
So they do not give any adjacency between these two triangles. Consequently, $T_0^i$ and $T_a^j$ are adjacent in  $G$ if and only if $i=j$.

We now use this to determine the link. The two surviving vertices in $T_0$ are adjacent in $G$, since $T_0$ is a triangle. The same holds for the two surviving vertices in $T_a$. Among the four cross-pairs formed by the surviving vertices, the two pairs with equal second coordinates are adjacent in $G$. The remaining two cross-pairs have different second coordinates
and are therefore non-adjacent in $G$.

By Lemma~\ref{lem:link-edge-criterion}, two vertices of $S$ form an edge of $\operatorname{lk}_{\Delta}(F)$ if and only if they are non-adjacent in $G$. Hence, the two remaining cross-pairs are precisely the
edges of the link. Since they have no common vertex, the link consists of two disconnected edges. Therefore
$\operatorname{lk}_{\Delta}(F)\cong K_2\sqcup K_2$.
\end{proof}

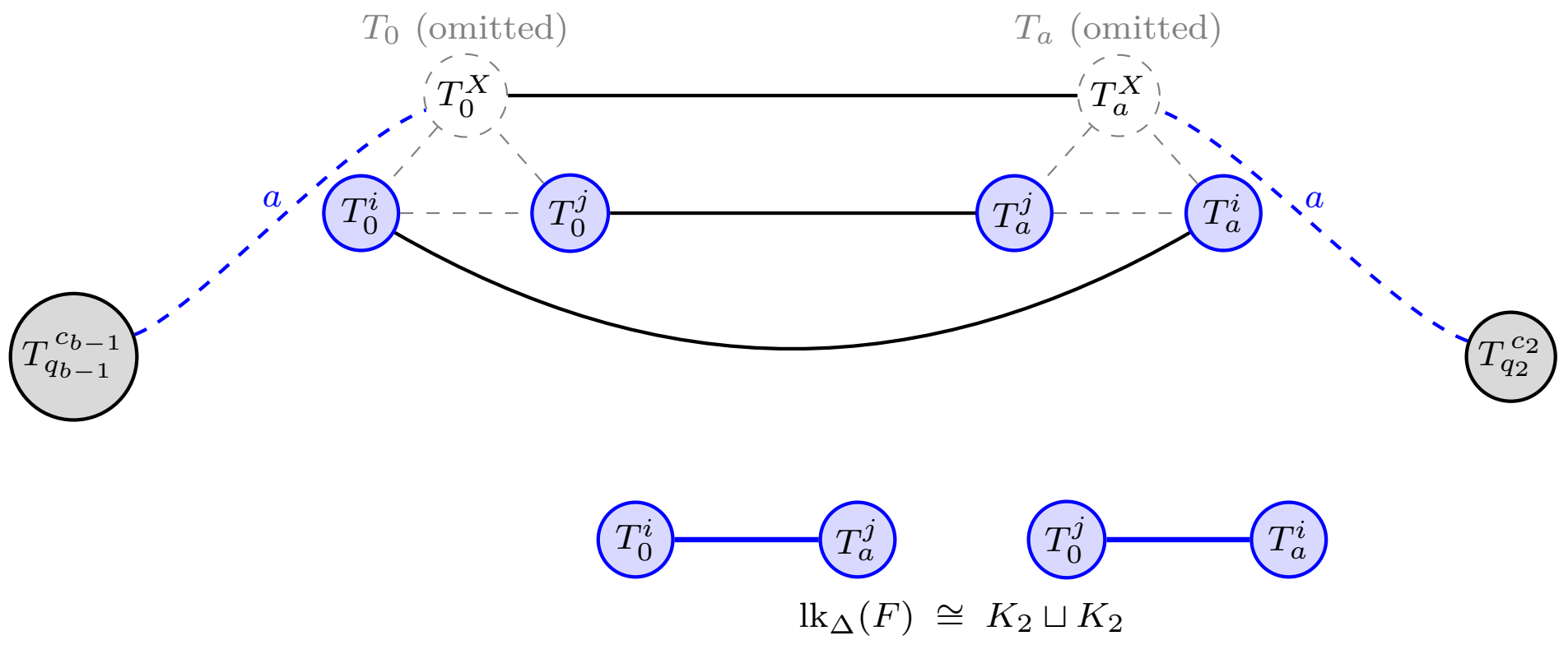
\begin{figure}[!ht]
\centering
\begin{tikzpicture}[scale=1.05,
  tri/.style={circle,draw,inner sep=1pt,minimum size=6mm,font=\scriptsize},
  sel/.style={tri,fill=black!15,thick},
  link/.style={tri,fill=blue!15,draw=blue,thick},
  omit/.style={tri,draw=gray,dashed}]

  \node[omit] (t0X) at (0,2)      {$T_0^{X}$};
  \node[link] (t0a) at (-0.8,1.1) {$T_0^{i}$};
  \node[link] (t0b) at (0.8,1.1)  {$T_0^{j}$};
  \draw[dashed,gray] (t0X)--(t0a)--(t0b)--(t0X);
  \node[font=\scriptsize,gray] at (0,2.5) {$T_0$ (omitted)};

  \node[omit] (taX) at (5,2)      {$T_a^{X}$};
  \node[link] (taj) at (4.2,1.1)  {$T_a^{j}$};
  \node[link] (tai) at (5.8,1.1)  {$T_a^{i}$};
  \draw[dashed,gray] (taX)--(taj)--(tai)--(taX);
  \node[font=\scriptsize,gray] at (5,2.5) {$T_a$ (omitted)};

  \draw[thick] (t0X) -- (taX);
  \draw[thick] (t0b) -- (taj);
  \draw[thick] (t0a) to[bend right=30] (tai);

  \node[sel] (tbm1) at (-3,0) {$T_{q_{b-1}}^{\,c_{b-1}}$};
  \draw[thick,dashed,blue] (tbm1) .. controls +(20:1) and +(200:1) .. (t0X)
       node[midway,above,font=\scriptsize,blue]{$a$};

  \node[sel] (t2) at (8,0) {$T_{q_2}^{\,c_2}$};
  \draw[thick,dashed,blue] (t2) .. controls +(160:1) and +(-20:1) .. (taX)
       node[midway,above,font=\scriptsize,blue]{$a$};

  \node[link] (l1) at (1.3,-1.4) {$T_0^{i}$};
  \node[link] (l2) at (3.0,-1.4) {$T_a^{j}$};
  \node[link] (l3) at (4.6,-1.4) {$T_0^{j}$};
  \node[link] (l4) at (6.3,-1.4) {$T_a^{i}$};
  \draw[very thick,blue] (l1)--(l2);
  \draw[very thick,blue] (l3)--(l4);
  \node[font=\scriptsize] at (3.8,-2.0)
       {$\operatorname{lk}_\Delta(F)\;\cong\;K_2\sqcup K_2$};

\end{tikzpicture}
\caption{Local structure near the two omitted triangles $T_0,T_a$. The face $F$ selects one vertex from every other triangle (shaded). The shown neighbours, whose layers $c_{b-1},c_2$ are chosen so that they are adjacent to $T_0^X,T_a^X$ respectively, remove these two vertices from the link. The solid edges $T_0^kT_a^k$ form the $a$-matching between $T_0$ and $T_a$, so the four surviving vertices (blue) are non-adjacent only across the pairs $T_0^iT_a^j$ and $T_0^jT_a^i$, giving $\operatorname{lk}_\Delta(F)\cong K_2\sqcup K_2$.}
\label{fig:general-disconnecting-face}
\end{figure}

\begin{theorem}\label{3b non CM}
Let $n=3b$, where $b\geq5$, and let $G=C_n(a,b)$ be the circulant graph satisfying $1\leq a<b\leq\left\lfloor\frac{n-1}{2}\right\rfloor$ and $\gcd(n,a,b)=1$.
Let $\Delta=\operatorname{Ind}(G)$ be the independence complex of $G$. Then $\Delta$ is not
Cohen--Macaulay over any field.
\end{theorem}

\begin{proof}
By the construction, the set $F=\{T_{q_k}^{c_k}:2\leq k\leq b-1\}$ is an independent set of $G$. Hence $F\in\Delta$. Then Lemma~\ref{Structure of the Link} $\operatorname{lk}_{\Delta}(F)\cong K_2\sqcup K_2$. Thus $\operatorname{lk}_{\Delta}(F)$ has exactly two connected components. Since it is the disjoint union of two edges, $\dim(\operatorname{lk}_{\Delta}(F))=1$. By Lemma~\ref{reduced-zero-homology}, for every field $K$, $\dim_K \widetilde H_0(\operatorname{lk}_{\Delta}(F);K)=2-1=1$.
Hence $\widetilde H_0(\operatorname{lk}_{\Delta}(F);K)\neq0$. However,$0<\dim(\operatorname{lk}_{\Delta}(F))=1$.
Therefore, for the face $F\in\Delta$ and $i=0$, we have $i<\dim(\operatorname{lk}_{\Delta}(F))$ but
$\widetilde H_i(\operatorname{lk}_{\Delta}(F);K)\neq0$. This contradicts Theorem~\ref{reisner}, which requires
$\widetilde H_i(\operatorname{lk}_{\Delta}(F);K)=0$ for every face $F\in\Delta$ and every $i<\dim(\operatorname{lk}_{\Delta}(F))$ when $\Delta$ is Cohen--Macaulay over $K$. Hence $\operatorname{Ind}(C_{3b}(a,b))$
is not Cohen--Macaulay over any field.
\end{proof}

\subsection{Setup and Structural Description for $n=3a$}
Let $n,a,b\in\mathbb{Z}_{>0}$,
$1\le a<b\le\left\lfloor\frac{n-1}{2}\right\rfloor$ and  $n=3a$, $a\geq4$, and $\gcd(n,a,b)=1$.

Let $\Delta=\operatorname{Ind}(G)$ denote the independence complex of
$G$. Put $r=b-a$. Since $a<b$, we have $r>0$. Moreover, as $n=3a$ and
$b<\frac{n}{2}=\frac{3a}{2}$, it follows that
$r=b-a<\frac{3a}{2}-a=\frac{a}{2}$. Thus, $0<r<\frac{a}{2}.$ Furthermore, $\gcd(a,r)=\gcd(a,b-a)=\gcd(n,a,b)=1.$

We next describe the structure of $G$ that is needed for the
construction of the face $F$. Since $n=3a$, every vertex $v\in\mathbb{Z}_n$
has a unique representation $v=t+ca$, where $0\leq t\leq a-1$ and $c\in\{0,1,2\}$.
We denote the vertex $t+ca$ by $T_t^c$. For each $t\in\{0,\ldots,a-1\}$, let
$T_t=\{T_t^0,T_t^1,T_t^2\}$. Thus, $T_t^0=t$, $T_t^1=t+a$, and $T_t^2=t+2a$, and the vertex set of $G$
is partitioned as $V(G)=T_0\cup T_1\cup\cdots\cup T_{a-1}$. The edges corresponding to the generator $a$ are considered, for every $t$,
$T_t^1-T_t^0=a$ and $T_t^2-T_t^1=a$, while $T_t^0-T_t^2\equiv-2a\equiv a\pmod{3a}$.
Since the connection set contains $\pm a$, these three vertices are pairwise
adjacent. Therefore, $G[T_t]\cong K_3$ for every $t$. Hence the vertices of $G$ are partitioned into $a$ triangles, $T_0,T_1,\ldots,T_{a-1}$. In particular, any independent set of $G$
can contain at most one vertex from each triangle $T_t$. This observation
will be used to exclude all edges arising from the generator $a$. It remains to understand how the edges corresponding to $b$ connect these
triangles. Since $b=a+r$, for a vertex $T_t^c=t+ca$ we have
$T_t^c+b=t+ca+a+r=t+r+(c+1)a$. The resulting expression depends on whether $t+r$ is less than $a$ or not.
If $t+r<a$, then $T_t^c+b=T_{t+r}^{\,c+1\pmod 3}$. On the other hand, if $t+r\geq a$, then
$t+r=(t+r-a)+a$, and therefore
$T_t^c+b=(t+r-a)+(c+2)a =T_{t+r-a}^{\,c+2\pmod 3}$. For convenience, define $e_t=
\begin{cases}
0, & t+r<a,\\
1, & t+r\geq a.
\end{cases}$

Then the two cases above can be written uniformly as \begin{equation}
\label{eq:b-action-setup}
T_t^c+b = T_{(t+r)\bmod a}^{\,c+1+e_t\pmod 3}.
\end{equation}

Thus, \eqref{eq:b-action-setup} describes precisely how an edge
corresponding to $b$ connects a vertex in one triangle to a vertex in another
triangle.\\

\begin{remark}
From $0<r<\frac{a}{2}$, we have $2r<a$. Recall that $e_t=0$ when
$t+r<a$ and $e_t=1$ when $t+r\geq a$. Hence $e_0=0$, since $r<a$, and
$e_r=0$, since $2r<a$. On the other hand, $(a-r)+r=a$, so
$e_{a-r}=1$. We shall also use the fact that $2r\notin\{0,r,a-r\}$. Indeed,
$2r\neq0$ and $2r\neq r$ because $r>0$. If $2r=a-r$, then $a=3r$.
Since $\gcd(a,r)=1$, this would imply $r=1$, and hence $a=3$,
contradicting $a\geq4$. Therefore $2r\neq a-r$. Thus,
$e_0=0$, $e_r=0$, $e_{a-r}=1$, and
$2r\notin\{0,r,a-r\}$.

\end{remark}

\begin{lemma} \label{F-independent}
Let $n=3a$, where $a\geq4$, and let
$G=C_n(a,b)$, where $b=a+r$ with
$0<r<\frac{a}{2}$ and $\gcd(a,r)=1$.
For $t\in\{0,\ldots,a-1\}$ and $c\in\{0,1,2\}$, write
$T_t^c=t+ca$, and put
$t^*=a-r$.
Define $F=
\{T_t^0:t\in\{1,\ldots,a-1\}\setminus\{r,t^*\}\}
\cup\{T_{t^*}^2\}.$ Then $F$ is an independent set of $G$. Consequently,
$F\in\Delta=\operatorname{Ind}(G)$.
Moreover,
$|F|=a-2$.
\end{lemma}

\begin{proof}
We prove that no two vertices of $F$ are adjacent in $G$. Recall that the edges of $G=C_n(a,b)$ arise from the differences
$\pm a$ and $\pm b$. Thus, it is sufficient to show that $F$ contains
neither an $a$-edge nor a $b$-edge. For each $t$, the three vertices
$T_t^0,T_t^1,T_t^2$ form a triangle, since any two of them differ by
$a$ modulo $n=3a$. Hence, the $a$-edges occur precisely within these
triangles. By the definition of $F$, at most one vertex is
chosen from each triangle. Therefore, no two vertices of $F$ can be
joined by an $a$-edge. We now consider all possible types of vertices occurring in $F$.

\medskip
\noindent\textbf{Case 1:} Two vertices of the form $T_t^0$.

Let $T_t^0\in F$. By \eqref{eq:b-action-setup},
$T_t^0+b=T_{(t+r)\bmod a}^{\,1+e_t\pmod 3}$.
Since $e_t\in\{0,1\}$, the second coordinate of this vertex is either
$1$ or $2$, and hence it is not of the form $T_s^0$. Therefore, a $b$-edge cannot join two vertices of $F$, both of the form
$T_t^0$.

\medskip
\noindent\textbf{Case 2:} A vertex of the form $T_t^0$ and the vertex $T_{t^*}^2$.

Suppose, for contradiction, that
$T_t^0+b=T_{t^*}^2$. Then the first coordinates must satisfy $t+r\equiv t^*=a-r\pmod a$. Hence
$t\equiv a-2r\pmod a $. Since $0<2r<a$, we obtain
$t=a-2r$. On the other hand, equality of the second coordinates requires $1+e_t\equiv2\pmod3,$ and therefore $e_t=1$. However, for $t=a-2r$, $t+r=(a-2r)+r=a-r<a$. Thus, $e_{a-2r}=0$, which contradicts the requirement $e_t=1$. Therefore no vertex of the form $T_t^0$ belonging to $F$ is
$b$-adjacent to $T_{t^*}^2$.

\medskip
\noindent\textbf{Case 3:} The image of the distinguished vertex $T_{t^*}^2$.

 Since $t^*=a-r$, we have $t^*+r=a$. Hence
$(t^*+r)\bmod a=0$ and $e_{t^*}=1$. Therefore, by
\eqref{eq:b-action-setup}, $T_{t^*}^2+b=T_0^1.$ But by the definition of $F$ in ,  $T_0^1\notin F$. Thus $T_{t^*}^2$ is not joined by a $b$-edge to any
other vertex of $F$.

The three cases above exhaust all possible pairs of vertices in $F$.
Hence $F$ contains neither an $a$-edge nor a $b$-edge. Therefore $F$ is
an independent set of $G$, and consequently
$F\in\Delta=\operatorname{Ind}(G)$. Finally, $|F|=(a-1-2)+1=a-2$.

\end{proof}

\begin{lemma}\label{lem:surviving-vertices}
Let $F$ be the independent set defined above. Then
$V(\operatorname{lk}_{\Delta}(F))
=\{T_0^0,T_0^2,T_r^0,T_r^1\}$.
\end{lemma}

\begin{proof}
We first show that every vertex of the link must belong to $T_0\cup T_r$. For every $t\notin\{0,r\}$, the set $F$ contains exactly one vertex from
the triangle $T_t$. Since each $T_t$ induces a triangle, every other
vertex of $T_t$ is adjacent to the vertex of $F$ lying in that triangle.
Consequently, no additional vertex from $T_t$ can be added to $F$ while
preserving independence. Hence no vertex of $T_t$, for
$t\notin\{0,r\}$, belongs to $\operatorname{lk}_{\Delta}(F)$. Therefore,
$V(\operatorname{lk}_{\Delta}(F))\subseteq T_0\cup T_r$.

Consider $T_0=\{T_0^0,T_0^1,T_0^2\}$.
Since $t^*=a-r$, we have $t^*+r=a$. Hence, by the structural relation
for the $b$-edges,
$T_{t^*}^2+b=T_0^1$.
Because $T_{t^*}^2\in F$, the vertex $T_0^1$ is adjacent to a vertex of
$F$. Therefore,
$T_0^1\notin V(\operatorname{lk}_{\Delta}(F))$. We now consider $T_0^0$ and $T_0^2$. The triangle $T_r$, which is the
image of $T_0$ under addition by $b$, contains no vertex of $F$.
Moreover, the only triangle whose $b$-edges enter $T_0$ is $T_{t^*}$,
and
$F\cap T_{t^*}=\{T_{t^*}^2\}$.
As observed above, $T_{t^*}^2$ is adjacent only to $T_0^1$ among the
vertices of $T_0$. Thus, neither $T_0^0$ nor $T_0^2$ is adjacent to any vertex of $F$. Since $F$ is independent, both vertices can be added to
$F$ while preserving independence. Hence
$T_0^0,T_0^2\in V(\operatorname{lk}_{\Delta}(F))$.
Consequently,
$V(\operatorname{lk}_{\Delta}(F))\cap T_0
=\{T_0^0,T_0^2\}$. Consider $T_r=\{T_r^0,T_r^1,T_r^2\}$.
Since $2r<a$, we have $e_r=0$. Therefore, the structural relation gives
$T_r^2+b=T_{2r}^0$. Moreover,
$2r\notin\{0,r,t^*\}$. Hence $T_{2r}$ is one of the triangles from which
$F$ contains the vertex $T_{2r}^0$. Thus $T_{2r}^0\in F$, and therefore
$T_r^2$ is adjacent to a vertex of $F$. Hence
$T_r^2\notin V(\operatorname{lk}_{\Delta}(F))$. It remains to consider $T_r^0$ and $T_r^1$. A $b$-edge entering $T_r$
must come from the triangle $T_0$, but $F\cap T_0=\varnothing$.
On the other hand, the $b$-edges starting from $T_r^0$ and $T_r^1$
go respectively to $T_{2r}^1$ and $T_{2r}^2$, whereas
$F\cap T_{2r}=\{T_{2r}^0\}$. Therefore, neither $T_r^0$ nor $T_r^1$ is
adjacent to any vertex of $F$. Hence, both can be added to $F$ while
preserving independence, and so
$T_r^0,T_r^1\in V(\operatorname{lk}_{\Delta}(F))$.
Consequently,
$V(\operatorname{lk}_{\Delta}(F))\cap T_r
=\{T_r^0,T_r^1\}$. Combining the two parts, we obtain
$V(\operatorname{lk}_{\Delta}(F))
=\{T_0^0,T_0^2,T_r^0,T_r^1\}$.
Thus, the vertices of the link are precisely the four vertices listed
above.
\end{proof}

\begin{lemma}\label{lem:edges-of-link}
Let $F$ be the independent set constructed above, then $\operatorname{lk}_{\Delta}(F)\cong K_2\sqcup K_2$.
\end{lemma}

\begin{proof}
By Lemma~\ref{lem:link-edge-criterion}, for any
$u,v\in V(\operatorname{lk}_{\Delta}(F))$, we have
$\{u,v\}\in\operatorname{lk}_{\Delta}(F)$ if and only if $u$ and $v$
are non-adjacent in $G$. By Lemma~\ref{lem:surviving-vertices},
$V(\operatorname{lk}_{\Delta}(F))
=\{T_0^0,T_0^2,T_r^0,T_r^1\}$. Since $n=3a$ and $b=a+r$, the connection set of $G$ is
$S=\{a,2a,a+r,2a-r\}$. Thus, two vertices of $G$ are adjacent
precisely when their difference modulo $n$ belongs to $S$. We now consider the six pairs of surviving vertices. The differences for
the pairs
$\{T_0^0,T_0^2\}$, $\{T_r^0,T_r^1\}$,
$\{T_0^0,T_r^1\}$, and $\{T_0^2,T_r^0\}$
are respectively $2a$, $a$, $a+r$, and $2a-r$. Each of these belongs to
$S$. Hence, these four pairs are adjacent in $G$, and therefore none
of them is an edge of $\operatorname{lk}_{\Delta}(F)$. For the pair $\{T_0^0,T_r^0\}$, the difference is $r$. Since
$0<r<a$, we have $r\neq a$ and $r\neq 2a$. Also, $r\neq a+r$, and
$r=2a-r$ would imply $r=a$. Therefore, $r\notin S$. Hence $T_0^0$ and
$T_r^0$ are non-adjacent in $G$, and consequently
$\{T_0^0,T_r^0\}\in\operatorname{lk}_{\Delta}(F)$. For the pair $\{T_0^2,T_r^1\}$, the difference is $a-r$. Since $r>0$,
we have $a-r\neq a$. Also, $a-r\neq 2a$ and $a-r\neq a+r$, while
$a-r=2a-r$ would imply $a=0$. Therefore, $a-r\notin S$. Thus
$T_0^2$ and $T_r^1$ are non-adjacent in $G$, and hence
$\{T_0^2,T_r^1\}\in\operatorname{lk}_{\Delta}(F)$. Therefore, the only edges of $\operatorname{lk}_{\Delta}(F)$ are
$\{T_0^0,T_r^0\}$ and $\{T_0^2,T_r^1\}$. Hence
\begin{equation}\label{eq:link-two-edges}
\operatorname{lk}_{\Delta}(F)\cong K_2\sqcup K_2.
\end{equation}
\end{proof}

\begin{theorem}\label{3a non CM}
Let $n=3a$, where $a\geq4$, and let $G=C_n(a,b)$ be the circulant
graph satisfying $1\leq a<b\leq\left\lfloor\frac{n-1}{2}\right\rfloor$ and $\gcd(n,a,b)=1$.
Let $\Delta=\operatorname{Ind}(G)$ be the independence complex
of $G$. Then $\Delta$ is not Cohen--Macaulay over any field.
\end{theorem}

\begin{proof}
By Lemma ~\ref{F-independent},
$F=\{T_t^0:t\in\{1,\ldots,a-1\}\setminus\{r,t^{*}\}\}
\cup\{T_{t^{*}}^2\}$. $F$ is an independent set of $G$.
Hence $F\in\Delta$. And then by Lemma ~\ref{lem:edges-of-link}
$\operatorname{lk}_{\Delta}(F)\cong K_2\sqcup K_2.$ Thus $\operatorname{lk}_{\Delta}(F)$ has exactly two connected
components. Moreover, since it consists of two disjoint edges,
$\dim\operatorname{lk}_{\Delta}(F)=1$. Therefore, by Lemma ~\ref{reduced-zero-homology},
for every field $K$, $\widetilde{H}_0(\operatorname{lk}_{\Delta}(F);K)
\cong K\neq0.$ Since
$0<\dim\operatorname{lk}_{\Delta}(F)=1$, we have found a face
$F\in\Delta$ and an integer $i=0$ satisfying
$i<\dim\operatorname{lk}_{\Delta}(F)$ for which
$\widetilde{H}_0(\operatorname{lk}_{\Delta}(F);K)\neq0$. Therefore, by Theorem ~\ref{reisner}
$\Delta=\operatorname{Ind}(G)$ is not Cohen--Macaulay over $K$.
\end{proof}
\begin{theorem}
Let $n,a,b\in\mathbb{Z}_{>0}$ satisfy
$1\leq a<b\leq\left\lfloor\frac{n-1}{2}\right\rfloor$ and
$\gcd(n,a,b)=1$ , and let
$G =C_n(a,b)$.  If $n=3a$ with $a\geq4$, or if $n=3b$ with $b\geq5$, then
$I(G)$ is not Cohen--Macaulay over any field.
\end{theorem}
\begin{proof}
    Follows from Theorems~\ref{3b non CM} and ~\ref{3a non CM}.
\end{proof}
\begin{corollary}
     Let $G=C_n(a,b)$ and for  $n=3a$ or $n=3b$,\ $G$ is never Cohen--Macaulay.

\end{corollary}
\begin{proof}
    In the verified range $n<19$ with $1\le a<b\le \left\lfloor \frac{n-1}{2}\right\rfloor, \gcd(n,a,b)=1,$
the circulant graphs $C_n(a,b)$ that are both well-covered and Cohen--Macaulay are exactly the following, $ C_5(1,2), C_7(1,2), C_7(1,3),C_7(2,3), C_8(1,2), C_8(2,3), C_{13}(1,5), C_{13}(2,3).$
From this we can observe that neither $n=3a$ nor $n=3b$ is Cohen--Macaulay.
\end{proof}

\section{Computational Verification}\label{appendix}

  \begin{center}
{\large\bfseries Well-Covered Circulant Graphs $C_n(a,b)$}\\[4pt]
{\small All $n < 19$,\quad $1 \le a < b \le \lfloor(n-1)/2\rfloor$,\quad
$\gcd(n,a,b)=1$}
\end{center}

\begin{longtable}{%
  >{\centering\arraybackslash}p{0.8cm}
  >{\centering\arraybackslash}p{2.5cm}
  >{\centering\arraybackslash}p{1.5cm}
  >{\centering\arraybackslash}p{1.2cm}
  >{\centering\arraybackslash}p{3.2cm}
  >{\centering\arraybackslash}p{2.0cm}
}

\toprule
$n$ &
Graph &
Well-Covered &
$\alpha(G)$ &
Maximal Ind.\ Set Sizes &
Cohen--Macaulay
\\
\midrule
\endfirsthead

\toprule
$n$ &
Graph &
Well-Covered &
$\alpha(G)$ &
Maximal Ind.\ Set Sizes &
Cohen--Macaulay
\\
\midrule
\endhead

\midrule
\multicolumn{6}{r}{\small\itshape Continued on next page\ldots}\\
\endfoot

\bottomrule
\endlastfoot

5 & $C_{5}(1,2)$ & Yes & $1$ & $\{1\}$ & Yes\\

6 & $C_{6}(1,2)$ & Yes & $2$ & $\{2\}$ & No\\

7 & $C_{7}(1,2)$ & Yes & $2$ & $\{2\}$ & Yes \\
7 & $C_{7}(1,3)$ & Yes & $2$ & $\{2\}$ & Yes \\
7 & $C_{7}(2,3)$ & Yes & $2$ & $\{2\}$ & Yes \\

8 & $C_{8}(1,2)$ & Yes & $2$ & $\{2\}$ & Yes \\
8 & $C_{8}(1,3)$ & Yes & $4$ & $\{4\}$ & No\\
8 & $C_{8}(2,3)$ & Yes & $2$ & $\{2\}$ & Yes \\

9 & $C_{9}(1,3)$ & Yes & $3$ & $\{3\}$ & No\\
9 & $C_{9}(2,3)$ & Yes & $3$ & $\{3\}$ & No\\
9 & $C_{9}(3,4)$ & Yes & $3$ & $\{3\}$ & No\\

10 & $C_{10}(1,4)$ & Yes & $4$ & $\{4\}$ & No\\
10 & $C_{10}(2,3)$ & Yes & $4$ & $\{4\}$ & No\\

11 & $C_{11}(1,2)$ & Yes & $3$ & $\{3\}$ & No\\
11 & $C_{11}(1,3)$ & Yes & $4$ & $\{4\}$ & No\\
11 & $C_{11}(1,4)$ & Yes & $4$ & $\{4\}$ & No\\
11 & $C_{11}(1,5)$ & Yes & $3$ & $\{3\}$ & No\\
11 & $C_{11}(2,3)$ & Yes & $4$ & $\{4\}$ & No\\
11 & $C_{11}(2,5)$ & Yes & $4$ & $\{4\}$ & No\\
11 & $C_{11}(3,4)$ & Yes & $3$ & $\{3\}$ & No\\
11 & $C_{11}(3,5)$ & Yes & $3$ & $\{3\}$ & No\\
11 & $C_{11}(4,5)$ & Yes & $4$ & $\{4\}$ & No\\

12 & $C_{12}(1,4)$ & Yes & $4$ & $\{4\}$ & No\\
12 & $C_{12}(3,4)$ & Yes & $4$ & $\{4\}$ & No\\
12 & $C_{12}(4,5)$ & Yes & $4$ & $\{4\}$ & No\\

13 & $C_{13}(1,3)$ & Yes & $5$ & $\{5\}$ & No\\
13 & $C_{13}(1,4)$ & Yes & $5$ & $\{5\}$ & No\\
13 & $C_{13}(1,5)$ & Yes & $4$ & $\{4\}$ & Yes \\
13 & $C_{13}(2,3)$ & Yes & $4$ & $\{4\}$ & Yes \\
13 & $C_{13}(2,5)$ & Yes & $5$ & $\{5\}$ & No\\
13 & $C_{13}(3,4)$ & Yes & $5$ & $\{5\}$ & No\\
13 & $C_{13}(5,6)$ & Yes & $5$ & $\{5\}$ & No\\

14 & $C_{14}(1,6)$ & Yes & $6$ & $\{6\}$ & No\\
14 & $C_{14}(2,5)$ & Yes & $6$ & $\{6\}$ & No\\
14 & $C_{14}(3,4)$ & Yes & $6$ & $\{6\}$ & No\\

15 & $C_{15}(1,5)$ & Yes & $5$ & $\{5\}$ & No\\
15 & $C_{15}(2,5)$ & Yes & $5$ & $\{5\}$ & No\\
15 & $C_{15}(3,5)$ & Yes & $5$ & $\{5\}$ & No\\
15 & $C_{15}(4,5)$ & Yes & $5$ & $\{5\}$ & No\\
15 & $C_{15}(5,6)$ & Yes & $5$ & $\{5\}$ & No\\
15 & $C_{15}(5,7)$ & Yes & $5$ & $\{5\}$ & No\\

18 & $C_{18}(1,6)$ & Yes & $6$ & $\{6\}$ & No\\
18 & $C_{18}(2,3)$ & Yes & $6$ & $\{6\}$ & No\\
18 & $C_{18}(3,4)$ & Yes & $6$ & $\{6\}$ & No\\
18 & $C_{18}(3,8)$ & Yes & $6$ & $\{6\}$ & No\\
18 & $C_{18}(5,6)$ & Yes & $6$ & $\{6\}$ & No\\
18 & $C_{18}(6,7)$ & Yes & $6$ & $\{6\}$ & No\\

\end{longtable}
\vspace{4pt}

{\footnotesize\textit{%
The independence complex $\Delta(G)$ was built in \textsf{Python} and
homology was computed via boundary-matrix rank over~$\mathbb{Q}$
using the Bron--Kerbosch algorithm~\cite{BronKerbosch1973}.
}}

\begin{remark}
The computational data for $n<19$ exhibit several genuine small-order
exceptions to the condition $n=3a$ or $n=3b$. Of the $47$ well-covered graphs, exactly $8$ are Cohen--Macaulay. Together with the theoretical results proved in the main text, these computations suggest
that $n=19$ may be a natural threshold for the eventual classification.
A complete verification for all $n\ge19$, however, is not established
here, since the parameter ranges not covered by Theorem~\ref{wellcovered} remain
open.
\end{remark}

\begin{theorem} \label{combined}
Let $n,a,b\in\mathbb{Z}_{>0}$ satisfy
$1\leq a<b\leq\left\lfloor\frac{n-1}{2}\right\rfloor$ and
$\gcd(n,a,b)=1$ with $n \geq 19$, and let
$G =C_n(a,b)$. For $c\in\{a,b\}$, set
$g_c=\gcd(c,n)$ and $m_c=n/g_c$. Throughout, when $b$ is the
reference generator, write
$g=g_b=\gcd(b,n)$ and $m=m_b=n/g$.

\begin{enumerate}
\item[(a)] If $n=3a$ or $n=3b$, then $G$ is well-covered and
$\alpha(G)=n/3$ for $n>6$.

Moreover, if $n=3a$ with $a\geq4$, or if $n=3b$ with $b\geq5$, then
$I(G)$ is not Cohen--Macaulay over any field.

\item[(b)] Suppose that $n\neq3a$ and $n\neq3b$. If either
\[
\begin{cases}
m\geq4\text{ is even}, & g\geq6,\\[4pt]
m\geq5\text{ is odd}, & g\geq4,
\end{cases}
\]
then $G$ is not well-covered. Consequently,
$I(G)$ is not Cohen--Macaulay over any field.
\end{enumerate}
\end{theorem}

\begin{proof}
\begin{enumerate}
  \item[(a)] Follows from Theorem~ \ref{thm:well-covered} that is G is well covered and follows from Theorem \ref{3b non CM} and Theorem \ref{3a non CM} that $I(G)$ is not Cohen--Macaulay over any field.
    \item[(b)] Follows from the Proposition \ref{prop:g5} and Proposition \ref{prop:g4} that G is not well-covered. Consequently, $I(G)$ is not Cohen--Macaulay over any field.
    \end{enumerate}
\end{proof}

\begin{remark}

The remaining cases are not completely classified by the
present results. More precisely, when $n\neq3a,3b$, the unresolved
well-coveredness cases include those for which
$m$ is even and $g\leq 6$, and those for which $m$ is odd and $g\leq3$.
In addition, for the exceptional well-covered cases $n=3a$ with $a\leq3$ and $n=3b$ with $b\leq 4$, the present Cohen--Macaulayness results do not determine whether
$I(G)$ is Cohen--Macaulay.
\end{remark}

\begin{theorem}\label{main}
    Let $d=gcd(n,a,b)$ , $G=C_n(a,b)$ and $H = C_\frac{n}{d} (\frac{a}{d}, \frac{b}{d})$ then,
    \begin{itemize}
        \item[(a)]  $G$ is connected if and only if $d=1$.
        \item[(b)] If $d>1$, then $G\cong d. H $ ,  $G$ is isomorphic to $d$ disjoint copies of $C_\frac{n}{d} (\frac{a}{d}, \frac{b}{d})$.
        \item[(c)]  $G$ is well-covered if and only if  $H$ is well covered.
        \item[(d)] $G$ is Cohen--Macaulay if and only if $H$ is Cohen--Macaulay.
        \end{itemize}
\end{theorem}
\begin{proof}
    Set $d=\gcd(n,a,b), a=da',b=db',n=dm,$
and define $H=C_m(a',b')=C_{n/d}(a/d,b/d).$
For integers $x,y$ and positive integers $d,m$, $ dx\equiv dy\pmod{dm} \iff x\equiv y\pmod m.$ Indeed, both sides say that $m\mid x-y$.

For $r=0,1,\dots,d-1$, define $ O_r=\{v\in\mathbb Z_n:v\equiv r\pmod d\}. $
Since $d\mid n$, these sets are well defined. They are pairwise disjoint and
$ \mathbb Z_n=O_0\cup O_1\cup\cdots\cup O_{d-1}. $ Each $O_r$ has exactly $m=n/d$ elements, namely $ r+kd, k=0,1,\dots,m-1. $ By basic observation if $u$ and $v$ are adjacent in $G$, we have  $ u-v\equiv 0\pmod d. $ Hence $ u\equiv v\pmod d.$
Therefore no edge of $G$ joins two different residue classes $O_r$ and $O_s$,
$\text{if }u\in O_r,\ v\in O_s,\ r\ne s,\text{ then }u,v\text{ are not adjacent.} $
\begin{itemize}

    \item[(a)] Assume $ d=\gcd(n,a,b)=1$ then there exist integers $x,y,z$ such that $ xa+yb+zn=1.$ Reducing modulo $n$, $ xa+yb\equiv 1\pmod n.$ Start at vertex $0$. Take $|x|$ steps, each adding $a$ if $x>0$, and subtracting $a$ if $x<0$. Then take $|y|$ steps, each adding $b$ if $y>0$, and subtracting $b$ if $y<0$. Every step moves along an edge of $G$, because consecutive vertices differ by $\pm a$ or $\pm b$. The final vertex is $ xa+yb\equiv 1\pmod n. $
Thus $0$ and $1$ lie in the same connected component of $G$. Now for every $j\in\mathbb Z_n$, the translation $\tau_j(v)=v+j$ is an automorphism of $G$, because adjacency depends only on the difference $u-v$. Applying $\tau_j$ to a path from $0$ to $1$ gives a path from $j$ to $j+1$. Hence $j$ and $j+1$ lie in the same connected component for every $j$. Therefore $ 0,1,2,\dots,n-1 $
all lie in one connected component. So $G$ is connected.

Assume $d>1$. The sets $O_0$ and $O_1$ are non-empty, since $0\in O_0, 1\in O_1.$
By the observation above, there is no edge between $O_0$ and $O_1$. Therefore, no path in $G$ can join $0$ to $1$. Hence $G$ is disconnected. Combining the two cases, $G$ is connected if and only if $d=1$.

\item[(b)] Assume $d>1$. By the observation above, there are no edges between different residue classes $ O_0,O_1,\dots,O_{d-1}.$ Therefore, $G$ is the disjoint union of the induced subgraphs
$ G[O_0],G[O_1],\dots,G[O_{d-1}].$ It remains to show that each $G[O_r]$ is isomorphic to $H$. Fix $r\in\{0,1,\dots,d-1\}$. Define $ \varphi_r:O_r\to\mathbb Z_m $ by $ \varphi_r(v)=\frac{\hat v-r}{d}\pmod m, $ where $\hat v$ is an integer representative of $v$. This map is well defined, changing $\hat v$ by a multiple of $n=dm$ changes $(\hat v-r)/d$ by a multiple of $m$. Its inverse is $ k\mapsto r+kd\pmod n, $ which is well defined because changing $k$ by a multiple of $m$ changes $kd$ by a multiple of $n=dm$. Hence $\varphi_r$ is a bijection. Now let
$ u=r+k_1d, v=r+k_2d $ be vertices in $O_r$. Then $ \varphi_r(u)=k_1,  \varphi_r(v)=k_2, $ and $ u-v\equiv (k_1-k_2)d\pmod{dm}.$ Thus $u$ and $v$ are adjacent in $G$ if and only if
$ (k_1-k_2)d\equiv \pm da'\pmod{dm} $ or $ (k_1-k_2)d\equiv \pm db'\pmod{dm}.$
By the congruence fact, this is equivalent to $ k_1-k_2\equiv \pm a'\pmod m $ or
$ k_1-k_2\equiv \pm b'\pmod m. $
But this is exactly the adjacency condition in $ H=C_m(a',b'). $
Therefore $\varphi_r$ is a graph isomorphism $ G[O_r]\cong H. $ Since this holds for every $r=0,1,\dots,d-1$, we obtain $ G\cong \underbrace{H\cup H\cup\cdots\cup H}_{d\text{ copies}}. $
Thus $G$ is isomorphic to $d$ disjoint copies of $H$.
\item[(c)]For $d=1$, the assertion is trivial because $G=H$. So assume $d>1$. we may write $ G=H_1\cup H_2\cup\cdots\cup H_d $ as a disjoint union of graphs with no edges between different $H_i$, where each $ H_i\cong H. $ For a vertex set $I\subseteq V(G)$, write $ I_i=I\cap V(H_i). $ Because there are no edges between different $H_i$, two vertices of $I$ can be adjacent only if they lie in the same $H_i$. Hence $I$ is independent in $G$ if and only if every $I_i$ is independent in $H_i$. Moreover, $I$ is a maximal independent set of $G$ if and only if every $I_i$ is a maximal independent set of $H_i$.

If every $I_i$ is maximal independent in $H_i$, then $I$ is independent. Let $v\notin I$, say $v\in V(H_i)$. Then $v\notin I_i$. By maximality of $I_i$ in $H_i$, the vertex $v$ has a neighbour in $I_i\subseteq I$. Hence $I$ is maximal in $G$. Conversely, if $I$ is maximal independent in $G$, then each $I_i$ is independent in $H_i$. Suppose some $I_i$ is not maximal in $H_i$. Then there exists a vertex $ v\in V(H_i)\setminus I_i $
such that $I_i\cup\{v\}$ is independent in $H_i$. Since $v$ has no neighbours outside $H_i$, the set $I\cup\{v\}$ is independent in $G$. This contradicts the maximality of $I$. Hence every $I_i$ is maximal in $H_i$.

Assume $H$ is well-covered. Let $s$ be the common size of all maximal independent sets of $H$. Let $I$ be any maximal independent set of $G$. By the description above, each $I_i$ is a maximal independent set of $H_i\cong H$. Therefore, $|I_i|=s$
for every $i$. Hence $ |I|=\sum_{i=1}^d |I_i|=ds. $ Thus every maximal independent set of $G$ has size $ds$. Therefore $G$ is well-covered.

Assume $G$ is well-covered. Let $A$ and $B$ be any two maximal independent sets of $H$. We must show $|A|=|B|.$ For each $i$, let $ \psi_i:V(H)\to V(H_i) $ be an isomorphism. Define
$ I^*=\psi_1(A)\cup\psi_2(A)\cup\cdots\cup\psi_d(A), $
and $
J^*=\psi_1(B)\cup\psi_2(A)\cup\cdots\cup\psi_d(A). $ Each $\psi_i(A)$ and $\psi_1(B)$ is a maximal independent set of the corresponding $H_i$, because $\psi_i$ is an isomorphism. By the description above, both $I^*$ and $J^*$ are maximal independent sets of $G$. Their sizes are
$ |I^*|=d|A|, $ and $ |J^*|=|B|+(d-1)|A|. $
Since $G$ is well-covered, $ |I^*|=|J^*|. $
Therefore $ d|A|=|B|+(d-1)|A|. $
Hence $ |A|=|B|. $ Since $A$ and $B$ were arbitrary maximal independent sets of $H$, all maximal independent sets of $H$ have the same size. Therefore $H$ is well-covered.

Combining all we have , $G$ is well-covered if and only if $H$ is well-covered.

\item[(d)] We may write $G=H_1\cup\cdots\cup H_d$ as a disjoint union of graphs with no edges between different $H_i$, where each $H_i\cong H$.
By the above observations, a set $I\subseteq V(G)$ is independent in $G$ if and only if $I_i:=I\cap V(H_i)$ is independent in $H_i$ for every $i$; since $I=I_1\cup\cdots\cup I_d$ is a disjoint union, this is equivalent to $I$ being a face of the join $\operatorname{Ind}(H_1)*\cdots*\operatorname{Ind}(H_d)$. Therefore
\begin{equation}\tag{1}
\operatorname{Ind}(G)=\operatorname{Ind}(H_1)*\cdots*\operatorname{Ind}(H_d),
\qquad
\operatorname{Ind}(H_i)\cong\operatorname{Ind}(H)
\quad\text{for }1\le i\le d.
\end{equation}

For $\sigma\in X$ and $\tau\in Y$, a face $\rho\cup\nu$ of $X*Y$, with $\rho\in X$ and $\nu\in Y$, belongs to $\operatorname{lk}_{X*Y}(\sigma\cup\tau)$ if and only if it is disjoint from $\sigma\cup\tau$ and $(\rho\cup\sigma)\cup(\nu\cup\tau)\in X*Y$. Because $X$ and $Y$ have disjoint vertex sets, this holds if and only if $\rho\in\operatorname{lk}_X\sigma$ and $\nu\in\operatorname{lk}_Y\tau$. Hence
\begin{equation}\tag{2}
\operatorname{lk}_{X*Y}(\sigma\cup\tau)
=
\operatorname{lk}_X\sigma*\operatorname{lk}_Y\tau.
\end{equation}
The faces of $X*Y$ correspond to pairs of faces of $X$ and $Y$, with a degree shift of one. This correspondence respects the boundary operators, and hence by Theorem~\ref{kunneth} we have
\begin{equation}\tag{3}
\widetilde H_j(X*Y;K)
\cong
\bigoplus_{p+q=j-1}
\widetilde H_p(X;K)\otimes_K \widetilde H_q(Y;K),
\text{ for all }j.
\end{equation}
Next we show that links of Cohen--Macaulay complexes are Cohen--Macaulay: if $\Delta$ is Cohen--Macaulay and $F\in\Delta$, then for $G\in\operatorname{lk}_\Delta F$ one has $\operatorname{lk}_{\operatorname{lk}_\Delta F}G=\operatorname{lk}_\Delta(F\cup G)$, so the conditions of \ref{reisner} for $\operatorname{lk}_\Delta F$ are among those for $\Delta$.

We claim that for complexes $X,Y$ on disjoint vertex sets,
\begin{equation}\tag{4}
X*Y\text{ is Cohen--Macaulay over }K
\iff
X\text{ and }Y\text{ are Cohen--Macaulay over }K.
\end{equation}
Now suppose $X$ and $Y$ are Cohen--Macaulay. Let $\sigma\cup\tau$ be a face of $X*Y$, with $\sigma\in X$ and $\tau\in Y$, and put $X'=\operatorname{lk}_X\sigma$, $Y'=\operatorname{lk}_Y\tau$. By (2), $\operatorname{lk}_{X*Y}(\sigma\cup\tau)=X'*Y'$, and $\dim(X'*Y')=\dim X'+\dim Y'+1$. Since links of Cohen--Macaulay complexes are Cohen--Macaulay, $X'$ and $Y'$ are Cohen--Macaulay, so $\widetilde H_p(X';K)=0$ for every $p<\dim X'$, and $\widetilde H_q(Y';K)=0$ for every $q<\dim Y'$. For $k<\dim X'+\dim Y'+1$, every summand in (3) with $p+q=k-1<\dim X'+\dim Y'$ has either $p<\dim X'$ or $q<\dim Y'$, hence vanishes. Therefore $\widetilde H_k\big(\operatorname{lk}_{X*Y}(\sigma\cup\tau);K\big)=0$. By \ref{reisner}, $X*Y$ is Cohen--Macaulay.

Conversely, suppose $X*Y$ is Cohen--Macaulay. Choose a facet $\tau$ of $Y$, which exists since $\emptyset\in Y$. Then $\operatorname{lk}_Y\tau=\{\emptyset\}$, because any $\rho\in\operatorname{lk}_Y\tau$ satisfies $\rho\cup\tau\in Y$ and maximality of $\tau$ forces $\rho=\emptyset$. By (2), $\operatorname{lk}_{X*Y}(\emptyset\cup\tau)=\operatorname{lk}_X\emptyset*\operatorname{lk}_Y\tau=X*\{\emptyset\}=X$, so $X$, being a link in a Cohen--Macaulay complex, is Cohen--Macaulay; by symmetry so is $Y$.

Now let $d\ge 1$. For $d=1$ there is nothing to prove, since $\operatorname{Ind}(G)\cong\operatorname{Ind}(H)$. For $d>1$, write $\operatorname{Ind}(G)=\operatorname{Ind}(H_1)*\Delta'$, where $\Delta'=\operatorname{Ind}(H_2)*\cdots*\operatorname{Ind}(H_d)$. By (4), $\operatorname{Ind}(G)$ is Cohen--Macaulay over $K$ if and only if both $\operatorname{Ind}(H_1)$ and $\Delta'$ are Cohen--Macaulay. By the induction hypothesis, $\Delta'$ is Cohen--Macaulay if and only if $\operatorname{Ind}(H)$ is Cohen--Macaulay, and $\operatorname{Ind}(H_1)\cong\operatorname{Ind}(H)$. Therefore, $\operatorname{Ind}(G)$ is Cohen--Macaulay over $K$ if and only if $\operatorname{Ind}(H)$ is Cohen--Macaulay over $K$; that is, $G$ is Cohen--Macaulay over $K$ if and only if $H$ is Cohen--Macaulay over $K$.
\end{itemize}
\end{proof}

\begin{remark}
    Theorem~\ref{main} shows that for the classification of well-covered circulant graphs and Cohen--Macaulay  $C_n(a,b)$, it is enough to treat the connected case
$\gcd(n,a,b)=1. $ Indeed, a disconnected graph $C_n(a,b)$ is well-covered and Cohen--Macaulay  exactly when the connected graph $C_{n/d}(a/d,b/d) $ is well-covered and Cohen--Macaulay respectively, whose $\gcd(n/d,a/d,b/d)=1$.
\end{remark}


\bibliographystyle{amsplain}
\bibliography{references}


\end{document}